\documentclass[10pt]{article}
\usepackage{ametsoc2col}
\RequirePackage{xcolor}
\usepackage{amsmath,amsfonts,amssymb,amsthm,cite,graphicx,subfig,verbatim,multicol,color,
}

\theoremstyle{definition}

\theoremstyle{remark}

\newcommand{\ZZ}{\mathbb{Z}}      

\newcommand{\diver}{\operatorname{div}}
\newcommand{\curl}{\operatorname{curl}}

\numberwithin{equation}{section}

\newcommand{\myabstract}
{We describe tornadogenesis and maintenance using the $3$-dimensional vortex gas model presented in \citet{chorin} and developed further in \citet{flandoli02}.   We suggest that high-energy, super-critical vortices in the sense of \citet{benjamin62}, that have been studied by \citet{fiedlerrotunno86}, have negative temperature in the sense of \citet{onsager49} play an important role in the model. We speculate that the formation of high-temperature vortices is related to the helicity inherited as they form or tilt into the vertical and their interaction with the surface and boundary layer. We also exploit the notion of self-similarity to justify power laws derived from observations of weak and strong tornadoes presented in \citet{cai,wurman00,wurman05}. Analysis of a Bryan Cloud Model (CM1) simulation of a tornadic supercell reveals scaling consistent with the observational studies.}

\date{\today}

\begin{document}

\title{\textbf{\large{Possible Implications of a Vortex Gas Model and Self-Similarity for Tornadogenesis and Maintenance}}}


    
    \author{\textsc{Doug Dokken,}
				\thanks{\textit{Corresponding author address:}
				Doug Dokken, Mail number  OSS201,
				Dept. of Math., Univ. of St. Thomas, 2115 Summit Ave., St. Paul, MN 55105.
				\newline{E-mail: dpdokken@stthomas.edu}}\quad\textsc{Kurt Scholz, and Mikhail M. Shvartsman}\\
\textit{\footnotesize{Mathematics Department, University of St. Thomas, Saint Paul, MN}}
\and
\centerline{\textsc{Pavel B\v{e}l\'{\i}k}}\\
\centerline{\textit{\footnotesize{Mathematics Department, Augsburg College, Minneapolis, MN}}}
\and
\centerline{\textsc{Corey K. Potvin}}\\
\textit{\footnotesize{Cooperative Institute for Mesoscale Studies, and NOAA/OAR/National Severe Storms Laboratory, Norman, OK, USA}}
\and
\centerline{\quad\textsc{Brittany Dahl}}\\
\centerline{\textit{\footnotesize{School of Meteorology, University of Oklahoma, Norman, OK, USA}}}
\and
\centerline{\quad\textsc{Amy McGovern}}\\
\centerline{\textit{\footnotesize{School of Computer Science, University of Oklahoma, Norman, OK, USA}}}
}



\ifthenelse{\boolean{dc}}
{
\twocolumn[
\begin{@twocolumnfalse}
\amstitle

\begin{center}
\begin{minipage}{13.0cm}
\begin{abstract}
	\myabstract
	\newline
	\begin{center}
		\rule{38mm}{0.2mm}
	\end{center}
\end{abstract}
\end{minipage}
\end{center}
\end{@twocolumnfalse}
]
}

\section{Introduction}
\label{sec:intro}

In a recent paper, \citet{cai} defined the pseudovorticity by $\zeta_\text{pv}=\frac{\Delta V}{L}$, where $\Delta V=|(V_r)_\text{max}-(V_r)_\text{min}|$ is the difference between the maximum and minimum radial velocity of the mesocyclone (rotating updraft) and $L$ is the distance between them. Mobile Doppler radar observations of past tornadic and nontornadic storms were filtered using a range of Cressman influence radii (lower-bounded by the native data resolution and upper-bounded by mesocyclone diameter) to obtain points $(\log(\varepsilon),\log(\zeta_\text{pv}))$, where $\varepsilon$ is the finest resolvable scale of the filtered radar data. Cai then calculated the regression line for each storm and found that steeper negative slopes are indicative of tornadic storms, and that the threshold slope for strong tornadoes in his sample was approximately $-1.6$. The regression lines Cai calculated strongly fit the data over scales between that of the mesocyclone core and that of the ``edge'' of the mesocyclonic tangential flow,  indicating a vorticity vs.~scale power law is valid over those scales. This suggests a power law for the decay of the vertical component of vorticity outside the solid-body mesocyclone core of the form
\begin{equation}
  \label{eq:cai_power_law}
  \zeta
  \propto
  r^b
\end{equation}
for some $b<0$, where $r$ is the radial distance from the axis of the vortex. Cai observed that it may be correct to interpret the exponent as a fractal dimension associated with the vortex. The vorticity power law, if valid, may extend to smaller (including tornadic) scales, but this could not be determined given the limited resolution of the radar observations used in Cai's study.

In a paper devoted to analyzing mobile radar data obtained from a tornado that occurred in Dimmit, Texas (June 2, 1995) and was rated F$2$--F$4$, \citet{wurman00} found $v\propto r^b$, i.e., the tangential winds outside the tornado core roughly fit the modified Rankine vortex model. They calculated the exponent $b$ and found it to vary from $-0.5$ to $-0.7$. From data obtained in the intercept of the Spencer, South Dakota F$4$ tornado (May 31, 1998), \citet{wurman05} calculated $b=-0.67$. Cai noticed that the (threshold) vorticity power law exponent for strongly tornadic mesocyclones in his sample differs from the velocity power law exponent calculated for the strong tornado in \citet{wurman00} by $1$, which is consistent with the vorticity being the curl of the velocity. The larger vorticity/velocity power law exponents found in strongly tornadic mesocyclones are consistent with the observation of \citet{trapp99} that ``parcels that nearly conserve angular momentum penetrate closer to the central axis of the tornadic mesocyclones, resulting in large tangential velocities.'' As noted by Cai, hurricanes exhibit a similar velocity power law and exponent outside their eyewall (see, e.g., \citet{miller67}). This suggests that roughly the same vorticity power law may apply over a range of atmospheric vortex scales. It is curious that the power laws obtained by \citet{wurman00,wurman05} are consistent with the results obtained earlier by \citet{lund-snow93} in a vortex simulator.

Power laws with a particular scaling exponent arise when a phenomenon ``repeats itself on changing scales'' (see, e.g., \citet{mandelbrot83,barnsley88,barenblatt96,barenblatt03}). This property is called self-similarity. We propose that to more completely understand strong atmospheric vortices requires further exploration of their self-similarity. Self-similarity, especially when exhibited across objects ranging in scale, can point to important properties of the underlying dynamics. We focus on the possible self-similarity of tornadoes in this paper. As will be shown, tornadoes appear to exhibit local self-similarity and local homogeneity suggesting fractal phenomena.

Some tracks left by high-energy vortices within a tornado are as narrow as $30$ cm.
Some of these paths appear to originate outside the tornado and intensify as they move into the tornado.  We identify these vortices as supercritical in the sense of \citet{fiedlerrotunno86}.
Their analysis the work of \citet{barcilon67, Burggraf77, benjamin62}, suggests to us that the super-critical vortex below a vortex breakdown 
has its volume and its length decrease as the  the energy of the super-critical vortex increases.  This would suggests that the entropy (randomness of the vortex) is decreasing when the energy is increased.
Hence the inverse temperature, which is the rate of change of the entropy with respect to the energy, of the vortex is negative.  This temperature is not related to the molecular temperature of the atmosphere.   These vortices would be barotropic, however their origin could very well be baroclinic.   Recent results suggest vorticity produced baroclinicly in the rear-flank downdraft and that then decends to the surface and is tilted into the vertical is linked to tornadogenesis. Once these vortices come in contact with the surface, and the stretching and surface friction related swirl (boundary layer effects) are in the appropriate ratio, then by analogy with the work of \citet{fiedlerrotunno86} the vortex would have negative temperature and the vortex would now be barotropic.

The paper is organized as follows. In Section \ref{sec:tornadogenesis} we describe the tornadogenesis problem. In Section \ref{sec:cm1}, motivated by the results of \citet{lund-snow93,wurman00,wurman05,cai}, we discuss a numerical experiment producing a time series of slopes of vorticity lines and its connection with previous results. In Section \ref{sec:self-similarity} we address self-similarity of tornadoes and mesocyclones and the way it might be manifested. We also give a heuristic argument supporting Cai's power law and its associated exponent for strong tornadoes. In Section \ref{sec:vortex_gases} we discuss the vortex gas theory of \citet{onsager49,chorin,flandoli02} in two and three dimensions and give arguments for its role in modeling tornadogenesis and tornado maintenance.  We discuss the influence of the boundary layer on the possible formation of negative temperature vortices.  We suggest that supercritical vortices have negative temperature.  One important point is that while in classical statistical mechanics of a molecular gas a large number of particles is assumed, for a three dimensional vortex gas one can study a single vortex and its properties use using an ensemble much like the ensembles used to model weather forecasts.

In Section \ref{sec:conclusions} we discuss suction vortices, give conclusions, and describe future work.  
\section{Tornadogenesis}
\label{sec:tornadogenesis}

The search to understand tornadogenesis invariably involves the question: ``Where does the vorticity in the tornado originate?'' The most penetrating studies of this question have lead to the study of two types of vorticity: barotropic vorticity and baroclinic vorticity. Barotropic vorticity is vorticity that exists in the ambient environment and is frozen in the fluid and stretched and advected by the fluid. Baroclinic vorticity is vorticity that is generated by density currents in the fluid and is stretched and advected by the fluid. Definitive discussions of the role of barotropic and baroclinic vorticity in tornadogenesis and the mathematical decomposition of vorticity into barotropic and baroclinic parts and its consequences are given in \citet{davies-jones82,davies-jones84,davies-jones96,davies-jones00,davies-jones06a,davies-jones06b,davies-jones08}. Both types of vorticity or combinations thereof have been suggested as the origins of the vorticity in tornadogenesis. Based on film footage of tornadoes showing sheets of precipitation spiraling into tornadoes, \citet{fujita73,fujita75} suggested a ``barotropic'' method called ``Fujita's recycling hypothesis.'' In this process the precipitation falling near the interface between the updraft and downdraft transported vorticity to the surface and into the tornado, and the precipitation-rich air was recycled into the thunderstorm updraft by the tornado.
The numerical model and experiment of \citet{davies-jones08} shows that tornadogenesis can take place by a purely barotropic process. \citet{markowski03} have also found through numerical studies that baroclinic vorticity can be important, if not dominant, in tornadogenesis in a recycling-type process as well. Additional evidence supporting the recycling hypothesis is the observation that downdrafts associated with tornadic storms are generally warmer than downdrafts associated with storms that were nontornadic; this would make the updraft more buoyant in the tornado-producing storms and less buoyant in storms that do not produce tornadoes (see \citet{markowski02a}).

Recent studies of radar data by \citet{markowski08} and numerical simulations by \citet{straka07} reveal arching vortex lines (or vortex tubes) in the rear flank of supercell storms. These vortex lines appear to be almost synonymous with supercell thunderstorms. If these vortex lines can be focused into a small region, as more and more vortex lines enter this region, viscous interactions between neighboring vortex lines are believed to lead to mergers. This can ultimately lead to the creation of a strong vortex.  Several theories have been given for the production of the arching vortex lines. In one theory vorticity lines (or rings) baroclinically generated around the rear flank downdraft are advected toward the updraft as they descend (see \citet{straka07,markowski08}). The downstream portions of the vorticity lines are subsequently lifted and stretched by the updraft while the upstream portions continue to descend, forming vortex arches. A second theory holds that horizontal shear across the rear flank gust front is the source of the vortex arches (see \citet{lee-wilhelmson97a,lee-wilhelmson97b,lee-wilhelmson00,trevorrow12}). This would create a vortex sheet that could be stretched and rolled up into a tornado vortex. It seems plausible that a combination of these two processes could be present, with a reconnection of vortex lines produced by the two different processes. However, there may be other vorticity sources as well. As discussed later, we hypothesize that the possible fractal dimension of $1.6$ found by Cai and implied in Wurman's observational studies comes from the interactions of vortices produced in these shear regions.

Observational analysis of videos of the tornadogenesis phase in large-diameter tornadoes forming under low-cloud bases (see, e.g., \citet{wadena10}) suggests that vortices (vortex lines) which enter the developing tornado make a partial revolution about the ambient tornado vortex before folding up and dissipating. The vortex gas model described in Section \ref{sec:vortex_gases} follows \citet{chorin} and
suggests this folding up is necessary to conserve energy as the vortex stretches and/or interacts with other vortices. In this process some energy is transmitted to much smaller scales, the so-called inertial range, beginning the Kolmogorov cascade to the viscous range and then dissipating as heat. However, as the vortex stretches before kinking up, much energy is transmitted to the ambient vortex as kinetic energy of the flow, and this increases the vorticity of the tornado. In this theory, as more and more vortices successively enter the developing tornado, the process repeats itself many times gradually increasing the vorticity of the tornado vortex. Subsequently, the vorticity of the ambient vortex is increased (assuming the vortex lines are produced uniformly), which can be further enhanced by stretching, eventually achieving quasi-equilibrium with its environment. During this process, energy is transferred from the smaller scales to the larger scales in an inverse energy cascade. As more vortices enter the ambient tornado vortex, larger vortices tend to form; the stronger vortices at the core and slightly weaker vortices wrapping around them. Visually, the resulting flow could manifest itself as multiple vortices or as a large single vortex. The stretching of the vortices that enter the tornado eventually leads to the dissipation of their vorticity in a Kolmogorov cascade. The process just described contrasts with the transition to a multiple-vortex configuration that can occur when the critical swirl ratio (a measure of tornado-scale helicity) is exceeded (see, e.g., \citet{church77}).

It has been shown by \citet{moffat69, moffat92} that twisting of subvortices about one another is measured by the helicity of the parent vortex. Helicity is the integral over physical space of the dot product of the velocity and the vorticity. It is thought that helicity of a flow inhibits the dissipation of energy and helps maintain the intensity of the flow (see \citet{lilly83,levich85,lilly86a,lilly86b}). Helicity and its slightly modified form have been used as a parameter to study its effect on supercell storms by \citet{lilly83,davies-jones84,lilly86b,davis-jones90,droegemeier93}. Lilly thought of a vortex as a coiled spring that unwinds as it stretches. If we think of the helicity as measuring how much the vortex is wound up, the stretching unwinds the spring and releases energy to the surrounding flow. This unwinding could manifest itself as vortex breakdown and/or the fractalization and kinking up that is predicted in the vortex gas theory (described in Section \ref{sec:vortex_gases}).

\citet{lee-wilhelmson97a,lee-wilhelmson97b,lee-wilhelmson00} studied non-super\-cell tornadogenesis due to vertical shear in the boundary layer. They considered a weak cold pool (outflow boundary) advancing from the west into an ambient flow from the south to the north. This led to a south-to-north oriented vortex sheet forming at the interface of the two flows. They noted first-generation vortices rolling up into stronger second-generation vortices. It seems plausible that if one did finer-grid simulations of the situation considered in the papers,
the first-generation vortices would have formed from roll-ups at smaller scales, resulting in a self-similar structure. The structure of the resulting vortex sheet resembled that of \citet{baker90}. The roll-up process has also been studied in \citet{snow78,rotunno84}. Snow concludes that ``the subsidiary vortices are integral parts of the overall flow pattern and should not be viewed as interacting independent vortices.'' Once  a 
vortex sheet roll-up has occurred, tornadogenesis can be induced by convection moving over and stretching one of the vortex sheet vortices. The first of three non-supercell tornadoes studied by \citet{roberts95} formed in this manner.

It is shown in \citet{roberts95} that with sufficient stretching, non-supercell tornado vortices can have EF3 strength, but virtually all violent (EF4--EF5) tornadoes occur within supercells. While this suggests important differences between tornadogenesis mechanisms in supercells versus non-supercells, it does not preclude the possibility that vortex sheet roll-up can also play an important role in supercell tornadogenesis. Radar analyses in \citet{dowell97,bluestein00a,bluestein00b,dowell02a,dowell02b,bluestein03a,bluestein03b} and observational analyses in \citet{brandes78,wilson86,wakimoto89,wakimoto96} have revealed vortices along the rear flank gust front of both tornadic and nontornadic supercell storms. The role of these vortices in tornado formation is a subject of current research. As with non-supercells, supercell tornadogenesis has been observed by \citet{wakimoto96} to sometimes be triggered when one of these vortices is stretched by an updraft (but this process alone probably cannot produce violent tornadoes). Consistent with observations, numerical supercell simulations by \citet{adlerman00} have shown vortices that appear to form along the edge of the rear flank gust front or a secondary gust front, and then roll up into a tornado vortex. The resulting vorticity distribution at one stage of tornadogenesis in \citet{adlerman00} resembles the two-dimensional vortex sheet roll-up modeled by \citet{chorin73} and \citet{krasny93}. The remarkable photo taken by Gene Moore (see Figure \ref{Vortex sheet roll-up along the rear flank gust front}) strongly implies vortex sheet roll-up. A sequence of vortices appears to be spiraling into a tornado as it crosses a lake. These feeder vortices have cross-sections too small to be resolved in all but the highest-resolution simulations currently achievable and appear to be very intense. Evidence of similar feeder vortices is also shown in Figure \ref{fig:suctionspots}, in which tracks left in corn fields appear to show vortices spiraling into the tornadoes and then dissipating as they stretch.

Subsequent to vortex sheet roll-up, vortex mergers may play a critical role. Observations of tornadogenesis near Bassett, NE by \citet{bluestein00a} suggested that a larger vortex ($\sim500$ m scale) and a smaller vortex ($\sim 100$--$200$ m scale), both of which were conjectured to arise from vortex sheet roll-up along the rear flank gust front, interacted so that the smaller vortex was absorbed by the larger vortex, possibly triggering tornadogenesis. The authors suggested that the origin of the vorticity was tilting of stream-wise vorticity along the rear flank gust front (see also \citet{dowell02b}). Tilting of stream-wise vorticity would result in vortices that have large helicity and are more resistant to dissipation due to stretching. In such a scenario, it is plausible that stretching and subsequent intensification of a larger vortex could draw other vortices through their mutual interaction, resulting in vortex mergers that further intensify the dominant vortex, and so on to tornadogenesis.

Further evidence that vortex sheet roll-up contributes to tornadogenesis in supercells is provided by \citet{chorin73}, which showed that vortex sheets consisting of cyclonically rotating vortices roll up into a cyclonic vortex. The rolled up vortex sheet resembled the hook echo region of a supercell thunderstorm. We believe the roll-up in \citet{chorin73} is representative of the roll-up in \citet{adlerman02}. When vortices of opposite sign were placed in the two halves of the vortex sheet segment, the sheet rolled up into a cyclonic--anticyclonic couplet resembling that often observed to straddle the hook echo and recently associated with arching vortex lines.

The process of tornadogenesis by roll-up of a vortex sheet undergoing stretching by the mesocyclone updraft occurs in the numerical simulation described in Section 3 and is illustrated in Figure \ref{vortex_image_time_series1.pdf}.  The time series for the maximum surface vorticity within the developing tornado can be found in Figure \ref{timeseries-ok-nssl}. Note how the vortex sheet vortices intensify as they stretch and approach the developing tornado. As the vortex sheet rolls up, the vortices transfer energy to the developing tornado vortex, increasing its maximum vorticity from $\sim0.1$ s$^{-1}$ to $\sim0.7$ s$^{-1}$ over the $30$-minute period shown. In Figure \ref{wurman-kosiba-fig7}, taken from \citet{wurman2013}, maximum gate-to-gate shear in an observed tornado exhibits marked oscillations superimposed on an upward trend. Given that the tornado did not display multiple-vortex behavior, we speculate that the shear oscillations and gradual increase in kinetic energy resulted from absorption of successive vortex sheet vortices, similar to that in the simulation. As such vortices are absorbed they would contribute not only energy but also helicity, which would decrease energy dissipation in the tornado (see \citet{andre77,lilly86a,lilly86b,yokoi93}).

\section{Vorticity lines computed from numerical simulation}
\label{sec:cm1}
A supercell thunderstorm simulation was investigated to help confirm the conclusions of \citet{cai} regarding the evolution of mesocyclone vorticity lines prior to and proceeding tornadogenesis. The supercell was simulated using the compressible mode of the non-hydrostatic Bryan Cloud Model 1 (CM1; \citet{bryanfritsch02}). The simulation proceeded on a $112.5$ km $\times$ $112.5$ km $\times$ $20.0$ km domain with horizontal grid spacing of $75$ m and vertical grid spacing increasing from $50$ m at the lowest layer to $750$ m at the highest layer. The large and small time steps were $1/4$ s and $1/16$ s, respectively. Typical of idealized storm simulations, a horizontally uniform analytical base state was used (see Figure \ref{skewt-ok-nssl}), terrain, surface fluxes, radiative transfer, and Coriolis acceleration were omitted, and radiative (free slip) lateral (vertical) boundary conditions were imposed. Microphysical processes were parametrized using the double-moment \citet{morrison05} scheme. The subgrid turbulence scheme was similar to \citet{deardorff80}. The simulated supercell exhibits features commonly observed in real supercells, including a hook echo reflectivity signature with a cyclonic--anticyclonic vorticity couplet (see Figure \ref{simulated_fields}).

In order to compute maximum vorticity at different length scales $\varepsilon$, the vorticity field valid on the $75$-m simulation grid was filtered using the \citet{cressman59} interpolation method with the cutoff radius set to $2\varepsilon$ (consistent with \citet{cai}). Vorticity lines were then computed near the low-level mesocyclone $\sim 500$ m above ground level (AGL) every $5$ minutes once a distinct low-level mesocyclone had formed (as discerned from visual inspect of the $75$-m vorticity field). As in \citet{cai}, vorticity lines were fit to $300$ m $\leq \varepsilon \leq 9600$ m. Tornadogenesis was considered to occur once the maximum axisymmetric tangential wind velocity, $V_T$, around the intensifying surface vortex associated with the low-level mesocyclone exceeded $20$ m s$^{-1}$. The $V_T$ was retrieved using the vortex detection and characterization technique of \citet{potvin13}.

As in \citet{cai}, the vorticity lines steepen prior to tornadogenesis (see Figures \ref{vorticitylines-ok-nssl} and \ref{timeseries-ok-nssl}), consistent with the concentration of vorticity from larger to smaller scales. Also as in \citet{cai}, a power law for vorticity appears to hold for scales exceeding that of the low-level mesocyclone core, but breaks down at smaller scales. \citet{cai} attributed this to smaller scales being more poorly resolved in his radar dataset. A similar effect occurs in our scenario: the effective model resolution of \citet{frehlich2008} artificially decreases the energy contained at scales approaching the grid spacing. In the absence of positive evidence that the vorticity power law indeed extends to tornadic and smaller scales, we can only offer this as a speculative explanation for the flattening of the vorticity lines at sub-mesocyclone scales. As discussed later, however, Cai's vorticity power law hypothesis (including its validity at tornadic scales) finds support in heuristic considerations of Kelvin--Helmholtz instability in vortex sheets.

\section{Self-Similarity }
\label{sec:self-similarity}
In this section we give possible ways mesocyclones and tornadoes might acquire self-similarity and so give rise to the hypothesized vorticity and velocity power laws discussed above.

Self-similarity can manifest itself in several ways in atmospheric flows. One such manifestation is scale-invariance of some characteristic of the flow, which may be demonstrated by the existence of a power law for the characteristic. Examples include the scenarios discussed above, where a power law for vorticity/pseudovorticity (\citet{cai}) or velocity (\citet{wurman00,wurman05}) are hypothesized.


Another manifestation is geometric self-similarity: features having similar shape occur at different scales. For example, the left image in Figure \ref{fig:hier-vort}, taken from \citet{church77}, illustrates a hierarchy of known vortex scales in tornadic supercells. This figure is strikingly similar to images in \citet[p.~168]{arnold99} and the Smale--Williams attractor cross-section in \citet{katok-hasselblatt95}, shown as the right image in Figure \ref{fig:hier-vort}. Videos of recent large tornadoes show subvortices of subvortices within tornadoes (see \citet{wadena10,elrenovideo13}). Though these sub-subvortices are transient and short lived, their existence suggests that near the surface the tornadic flow is approximating a vertically periodic flow similar to the Smale--Williams attractor. Recently, \citet{kuznetsov11} produced a dynamical system with a Smale--Williams attractor. The system was forced by periodic pulses. We will revisit the idea of periodic pulses in the conclusions.

Geometric self-similarity is occasionally seen in high-resolution numerical simulations of tornadic supercells (\citet{nova04,adlerman00}) and also in Doppler radar and reflectivity observations (\citet{bluestein00b,nova04}); see Figure \ref{fig:tornado_fractal} for an example. High-quality video recordings of some recent tornadoes depict mini-suction vortices (subvortices of suction vortices), confirming the smallest scale of the hierarchy in Figure \ref{fig:hier-vort}.

In a related work, \citet{belik14} revisit the ``swirling vortex'' model of \citet{serrin} and investigate solutions to the Navier--Stokes and Euler equation

$v= r^b$, where $b$ is not necessarily equal to $-1$. The streamlines of the modeled vortices exhibit self-similarity, i.e., both the power law and the geometric manifestations of self-similarity are addressed in this work.

\citet{chorin} in his study of turbulent flows found quantities with fractal dimensions. In numerical experiments he found the fractal dimensions of the axes of vortices he studied to be related to the ``temperature'' of the vortex. ``Hot'' negative-temperature vortices had a smooth axis, while at temperatures of positive or negative infinity (Kolmogorov cascade) the vortex had a fractal axis (and cross-section). We hypothesize that high-energy vortices entering the tornado acquire fractal axes upon being stretched and kinked up (transition from negative-temperature vortices to infinite-temperature vortices). One would expect a mixture of fractal dimensions for these axes in the turbulent region surrounding the solid body tornado core.

In their study of the effect of rotation and helicity on self-similarity, \citet{pouquet10} comment that ``when comparing numerical simulations, it was found that two runs at similar Rossby number and at similar times (albeit at different Reynolds number) display self-similar behavior or decreased intermittency depending on whether the flow had helicity or not.'' That tornadoes form in helical environments may largely account for the degree of self-similarity that is often observed in them (e.g., the presence of suction vortices), and suggests self-similarity may extend to smaller scales than currently known.  We propose that such self-similarity can arise within persistent vortex sheets along the rear flank downdraft gust fronts of tornadic supercells. In the proposed scenario, a sequence of vortex roll-ups occurs, with each new generation of vortices forming from previous-generation vortices wrapping around each other, ultimately resulting in vortices with roughly fractal cross-sections (geometric self-similarity).

We now give a heuristic argument to support Cai's power law for the vertical vorticity of strong tornadoes, i.e., $\zeta=\mathcal{O}(\epsilon^{\sim(-1.6)})$. From Kelvin's circulation theorem, the product of vorticity $\zeta$ and the cross-sectional area of a vortex tube is constant for Eulerian barotropic flows. Hence, $\zeta=C/A$, where $A$ is the cross-sectional area of the vortex. A numerical study of Kelvin--Helmholtz instability in \citet{baker90} identifies a relationship between the thickness of the vortex sheet, $h$, and the cross-sectional area of the vortices, $A$. The result is that $A$ scales like $\mathcal{O}(h^{1.55})$ as $h\to0$. To the degree that a tornado is formed from vortex sheet vortices, and that the hypothesized vorticity power law (and self-similarity) extends to tornadic and smaller scales, $\zeta$ would then scale as $\mathcal{O}(h^{-1.55})$. If the vortices stretch within the updraft, their cross-sections decrease and their vorticity increases, causing the slope of the vorticity line to decrease. Under this scenario, $-1.55$ would provide an upper bound for the vorticity line slope of a tornado forming by vortex sheet roll-up. That tornadoes may derive much of their energy from vortex sheet vortices is made plausible by the fact that energy can cascade from smaller to larger scales in two-dimensional flows (\citet{chorin}). Vortex sheet vortices and tornadoes (at least near the surface) are approximately two-dimensional.

\section{Vortex gases}
\label{sec:vortex_gases}
In this section we give an overview of a vortex gas theory for two dimensions and a vortex gas model for three dimensions. We introduce a notion of entropy and temperature that is different from the usual notions of entropy and temperature of gases of molecules. We use the theory and models described here to address the question of tornadogenesis and maintenance. The interaction of large numbers of vortices in two- and three-dimensional space has been studied by modeling the vortices as part of a vortex gas. This theory has its origins in the $19$th century in the works of \citet{helmholtz58} and \citet{kelvin69}. The theory is the analogue of the classical statistical mechanics of gases, which attempts to explain the macroscopic behavior of gases by using the statistics of  modeled microscopic   behavior of molecules. In the vortex gas case the molecules are replaced by vortices. These could also be arching vortex lines (tubes). Just as in the case of gases, our development includes specialized entropy and temperature.  It is important to note that one does not need an inordinately large number of vortices to use this theory.  \citet{onsager49} first suggested the notion of temperature for vortex gases and formulated a two-dimensional theory. For a discussion of these ideas see \citet{chorin93,chorin}. We follow the development of \citet{chorin-marsden93,chorin,newton01}, being guided by the work and energy balance analysis of \citet{lilly86a} and finding some analogues in turbulence theory that one could apply to the development of rotation in tornadoes. In classical thermodynamics a body in contact with a heat bath will heat up until it achieves equilibrium with its environment. By analogy, a vortex will heat up if it interacts with ``hotter'' vortices. As this process repeats itself many times, the vortex eventually achieves quasi-equilibrium with its environment.  The frequency and intensity of these the "hotter" vortices vortices will determine the intensity of the tornado.  These ``hot'' vortices will have negative temperature (in the sense of vortex gas theory) and will increase the energy of a developing vortex. In the case where the vortex is a tornado or pre-tornado, this gives us a mechanism for understanding some aspects of tornadogenesis and maintenance. A critique of the notion of negative temperature in the two-dimensional vortex gas theory has been given by \citet{frohlich-ruelle82} and \citet{miller92}.  One does not need a large number of vortices to use the approximations.

The modeling of vortices in three dimensions has been carried out by \citet{chorin-akao91} using the Ising model. It uses this simplified approach to the vortex gas to study the relationship of stretching and temperature of the vortex and other quantities associated to a vortex gas. In this approach the vortices appear as either horizontal or vertical segments joining adjacent points in a three-dimensional lattice. The lattice is formed from the points in $\ZZ^3$,  the three-dimensional space with integer coordinates. As time advances, the vortex configuration is allowed to change without allowing the vortices to self-intersect. The future configurations of the vortices are then studied using a Monte Carlo Markov chain algorithm.   The integer lattice could be replaced by a lattice with smaller grid spacing for finer resolution and a closer vortex approximation.  Chorin uses ensembles to approximate the behavior of a vortex, similar to the ensembles of models used in weather modeling.  Motivated by the ideas of Chorin, \citet{flandoli02} developed a stochastic theory for vortex filaments in three dimensions which included fractal cross sections.  This theory was not restricted to vortices on a lattice, but more general vortices in three dimensions.  Based on the notion of capacity, cross sections of the vortices had to be fractal for the vortices to have finite energy.    The capacity of the vortex is related to its energy and as the fractal dimension of the cross section of the vortex increases, the capacity increases as well. 

We do not claim tornadic vortices are exactly modeled by an Ising model.  However the model can give insight into the qualitative behavior of vortices.  While Chorin's Ising model may seem unrelated to the case of tornadic vortices, the ideas of Flandoli and Gubinelli support ideas in this paper and Chorin's experiments give qualitative support to the paper.

In the last part of this section we take into account the cross-sections of vortices and combine the results of \citet{wurman00,wurman05,cai} with an argument of \citet{chorin} and use it to obtain information about the possible source of increase in vorticity at tornado scales. 

\subsection{The two-dimensional vortex gas theory}
The Euler equation for incompressible fluid flow is 
\begin{equation*}
  \frac{D{\bf V}}{Dt}
  =
  \frac{\partial{\bf V}}{\partial t}
  +  ({\bf V}\cdot\nabla){\bf V}
  =
  -\nabla p+{\bf f},
\end{equation*}
where $\bf V$ is the velocity, $p$ is the pressure, and $\bf f$ an external body force.
To obtain the equation for vorticity, ${\boldsymbol\omega}=\curl{\bf V}=\nabla\times{\bf V}$, we take the curl of the above equation and obtain
\begin{equation*}
  \frac{D\boldsymbol\omega}{Dt}
  =
  (\boldsymbol\omega\cdot\nabla){\bf V}.
\end{equation*}
It is possible to extract the $z$ (vertical) component of vorticity, $\zeta$, from the above equation to obtain (\citet{klemp87})
\begin{equation}\label{vert-vort}
  \frac{\partial\zeta}{\partial t}
  =
  -{{\bf V}\cdot\nabla\zeta}
  +
  \zeta\,\frac{\partial w}{\partial z}
  +
  {\boldsymbol\omega}_H\cdot{\nabla_Hw},
\end{equation}
where $w$ is the vertical component of the velocity, $\boldsymbol\omega_H$ is the horizontal component of the vorticity, and $\nabla_H$ is the horizontal gradient. The three terms on the right-hand side of the above equation represent the change in vertical vorticity due to the advection, stretching and tilting of vorticity, respectively.

Vortices that form in strongly sheared environments have a two-dimensional structure before they are stretched. Recently, \citet{nolan12} has shown that a significant amount of the perturbation energy in tornadoes is due to stretching by the updraft near the surface. To model the flow behavior of an intense vortex at the surface, we assume the flow is essentially two-dimensional, although this does not fully capture three-dimensional behavior of the vortex. Under the influence of strong rotation, turbulent flow becomes anisotropic and the flow tends to become two-dimensional but never quite reaches that state (\citet{pouquet10}). 
Comparing the tracks left by suction spots in tornadoes (Figure \ref{fig:suctionspots}) and plots of interacting two-dimensional vortices (Figure \ref{fig:f-pics}) suggests that there is a connection between the two. The tracks suggest that the vortices behave like two-dimensional vortices, then dissipate due to stretching. This phenomenon can be observed in videos of intense tornadoes.  The existence of such vortices and their behavior may be related to the vortices studied in papers by \citet{fiedler94, fiedlerrotunno86, fiedler97, 
xialewellens03, lewellens07}.  \citet{fiedlerrotunno86} show that vortices in contact with ground that undergo stretching and an appropriate addition of swirl at the surface can become supercritical.
The behavior of the vortices is related to the swirl ratio.  The radius of the vortex is related to the thickness of the boundary layer and one thinks of the vortex erupting form the surface as an extension of the boundary layer.  If the stretching is dominant then the vortex will stretch and eventually dissipate.  \citet{fiedlerrotunno86} show if the swirl is dominant the supercritical vortex will narrow as the swirl increases and will become shorter.  As this occurs the vortex also intensifies  until a critical threshold is reached and the vortex undergoes a vortex breakdown.  
The intensification of the vortex, as swirl is added, is consistent with an increase in energy, while the decrease in length and radius is consistent with a decrease in entropy.  This suggests that the vortices have negative inverse temperature in the vortex gas sense, $dS/dE<0$, that is, when energy is added to the vortex, its entropy (randomness of the vortex) decreases.  If stretching dominates, the supercritical vortex would decrease in temperature and lead to its dissipation.  If the swirl is dominant,
then the vortex reaches a maximum temperature consistent with ideas of Joyce and Montgomery (\citet{chorin}).  The vortex gas theory does not say what will happen with an increase in energy beyond this value, however based on the analogy above, vortex breakdown appears to be the outcome.

We proceed to develop the two-dimensional theory to model the two-dimensional behavior. We then develop the three-dimensional model to fully understand dissipation. The two-dimensional vorticity equation for incompressible fluid flow is
\begin{equation*}
  \frac{D\boldsymbol\omega}{Dt}
  =
  0.
\end{equation*}
Writing the velocity in the component form, ${\bf V}=(u,v)$, the incompressibility condition, $\diver{\bf V}=0$, together with the assumption that the underlying domain is simply connected (has no holes) implies that there exists a stream function, $\psi(x,y)$, such that
\begin{equation}
  u
  =
  \frac{\partial\psi}{\partial y},
  \qquad
  v
  =
  -\frac{\partial\psi}{\partial x}
  \label{eq:stream_velocity}
\end{equation}
and
\begin{equation}
  -\Delta\psi
  =
  \zeta.
  \label{eq:poisson_for_psi}
\end{equation}
Assume that vorticity is concentrated at discrete points ${\bf x}_i=(x_i,y_i)$ for $i=1,\dots,n$, each with circulation $\Gamma_i$, so that
\begin{equation*}
  \zeta({\bf x})
  =
  \sum_{i=1}^n \Gamma_i \delta({\bf x-x}_i),
\end{equation*}
where $\delta$ denotes the Dirac delta function and ${\bf x}=(x,y)$. The solution to \eqref{eq:poisson_for_psi} is given by
\begin{equation*}
  \psi({\bf x})
  =
  -\sum_{i=1}^n\frac{\Gamma_i}{2\pi} \log{\|{\bf x-x}_i\|},
\end{equation*}
and, using \eqref{eq:stream_velocity}, the velocity field induced by the $j$th vortex is,
\begin{equation*}
  V_j({\bf x})
  =
  \frac{\Gamma_j}{2\pi r^2}\left(y_j-y,-(x_j-x)\right),
\end{equation*}
where $r=\|{\bf x-x}_j\|$.

If one assumes that each of the vortices moves under the influence of the combined velocity field of the remaining vortices, then
\begin{equation*}
  \frac{d{\bf x}_i}{dt}
  =
  \sum_{\substack{j=1\\j\ne i}}^nV_j({\bf x}_i)
  =
  \sum_{\substack{j\ne i}}^nV_j({\bf x}_i),
\end{equation*}
or, in the component form,
\begin{align*}
  \frac{dx_i}{dt}
  &=
  \frac{1}{2\pi}\sum_{j\ne i}\frac{\Gamma_j({y_j-y_i})}{r_{ij}^2}\\
  \frac{dy_i}{dt}
  &=
  -\frac{1}{2\pi}\sum_{j\ne i}\frac{\Gamma_j({x_j-x_i})}{r_{ij}^2},
\end{align*}
where $r_{ij}=\|{\bf x}_i-{\bf x}_j\|$. These equations form a Hamiltonian system that has rigorous connections with the Euler equation (\citet{chorin-marsden93,marchioro94}).
The corresponding Hamiltonian is
\begin{equation*}
  H
  =
  -\frac{1}{4\pi}\sum_{\substack{1\le i,j\le n\\i\ne j}}\Gamma_i\Gamma_j\log{\|{\bf x}_i-{\bf x}_j\|}.
\end{equation*}
It is easy to check that the Hamiltonian is conserved, i.e.,
\begin{equation*}
  \frac{dH}{dt}
  =
  0,
\end{equation*}
which, in particular, implies that if all the circulations are of the same sign, then the vortices cannot merge in finite time.
Other conserved quantities are the total vorticity, $\Gamma$, the center of vorticity, $\bf M$, and the moment of inertia, $I$, given by 
\begin{equation*}
  \Gamma
  =
  \sum\Gamma_i,
  \qquad
  {\bf M}
  =
  \frac{\sum\Gamma_i\bf{x}_i}{\Gamma},
  \qquad
  I
  =
  \sum\Gamma_i\|{\bf x}_i-{\bf M}\|^2,
\end{equation*}
where all the sums are for $i=1,\dots,n$.

Using this representation, one can model the behavior of vortex configurations in the plane (\citet{marchioro94,newton01,lim-nebus07,chorin73}). For example, a pair of vortices of equal
circulations will move about the midpoint of the segment joining them.
For a line of vortices of equal circulations, the vortices stay in a line.
If one has a half-line of vortices located at integer points on the $x$-axis, 
the vortex half-line rolls up into a spiral.

We next give a brief exposition of a two-dimensional theory of vortex gases by proceeding in analogy with the development of the Boltzmann distribution in the theory of statistical mechanics in three dimensions. The particles are replaced by vortices and the assumptions on the distribution of vortices in a region in two-dimensional space is used to define a distribution in the corresponding phase space. The entropy of the distribution is defined. The discussion that follows is general and can be used in the context of both two- and three-dimensional flows.

Consider a vortex system, with $k$ possible energy and moments of inertia levels $(E_j,I_j)$, $j=1,\ldots,k$. Let $0\le p_j\le1$ represent the probabilities that the system is in the energy state $E_j$ with the moment of inertia $I_j$. We assume that the average energy of the ensemble, $\langle E\rangle$, and the average moment of inertia, $\langle I\rangle$, are fixed. Then
\begin{equation}
  \label{constr1}
  \sum_{j=1}^k p_j
  =
  1,
\end{equation}
\begin{equation}
  \label{constr2}
  \sum_{j=1}^k p_j E_j
  =
  \langle E\rangle,
\end{equation}
\begin{equation}
  \label{constr3}
  \sum_{j=1}^k p_j I_j
  =
  \langle I\rangle.
\end{equation}
We define the entropy of the ensemble corresponding to the macrostate $(\langle E\rangle,\langle I\rangle)$, up to an additive constant, as
\begin{equation}
  \label{eq:entropy}
  S
  =
  -\sum_{j=1}^k p_j\log{p_j}.
\end{equation}
To maximize the entropy, we consider the Lagrangian
\begin{equation}
  \label{unconstrained}
  \begin{split}
    & L(p_1,\ldots,p_k)
      =
      -\sum_{j=1}^k p_j\log{p_j}
      -\alpha\left(\sum_{j=1}^k p_j-1\right)\\
    &\quad
      -\beta\left(\sum_{j=1}^k p_j E_j-\langle E\rangle\right)
      -\gamma\left(\sum_{j=1}^k p_j I_j-\langle I\rangle\right),
  \end{split}
\end{equation}
where $\alpha$, $\beta$ and $\gamma$ are Lagrange multipliers. Differentiating \eqref{unconstrained} with respect to each of the $p_j$ and setting these partial derivatives equal to zero, we obtain
\begin{equation*}
  -\log{p_j}-1-\alpha-\beta E_j-\gamma I_j
  =
  0,
  \qquad
  j=1,\dots,k,
\end{equation*}
while the derivatives with respect to the Lagrange multipliers return the constraints \eqref{constr1}--\eqref{constr3}. This results in
\begin{equation}
  \label{pj}
  p_j
  =
  e^{-1-\alpha-\beta E_j-\gamma I_j},
  \qquad
  j=1,\dots,k.
\end{equation}
From \eqref{constr1} and \eqref{pj} we now have
\begin{equation}
  \label{eq:partition_function}
  e^{1+\alpha}
  =
  \sum_{j=1}^k e^{-\beta E_j-\gamma I_j}
  \equiv
  Z,
\end{equation}
where $Z$ is called a partition function. It follows that
\begin{equation*}
  p_j
  =
  \frac{e^{-\beta E_j-\gamma I_j}}{Z},
\end{equation*}
and from \eqref{constr2} and \eqref{constr3} we have
\begin{equation*}
 \langle E\rangle
 =
 \sum_{j=1}^k E_j \frac{e^{-\beta E_j-\gamma I_j}}{Z},
 \quad
 \langle I\rangle
 =
 \sum_{j=1}^k I_j \frac{e^{-\beta E_j-\gamma I_j}}{Z}.
 \end{equation*}
Consequently, using the definition of the partition function \eqref{eq:partition_function}, we obtain
\begin{equation}
  \label{eq:partials_of_Z}
  \begin{split}
    -\frac{\partial\log{Z}}{\partial\beta}
    &=
    \frac{\sum_{j=1}^k E_j e^{-\beta E_j-\gamma I_j}}{Z}
    =
    \langle E\rangle,\\
    -\frac {\partial\log{Z}}{\partial\gamma}
    &=
    \frac{\sum_{j=1}^k I_j e^{-\beta E_j-\gamma I_j}}{Z}
    =
    \langle I\rangle.
  \end{split}
\end{equation}
The expression \eqref{eq:entropy} for the entropy can be written as
\begin{equation*}
  S
  =
  \sum_{j=1}^k p_j(\beta E_j+\gamma I_j+\log{Z})
  =
  \beta\langle E\rangle
  +
  \gamma\langle I\rangle
  +
  \log{Z},
\end{equation*}
and differentiating it with respect to $\langle E\rangle$ and using \eqref{eq:partials_of_Z} gives
\begin{equation}
  \label{beq1T}
  \begin{split}
    \frac{\partial S}{\partial\langle E\rangle}
    &=
    \frac{\partial\beta}{\partial\langle E\rangle}\langle E\rangle
    +
    \beta
    +
    \frac{\partial\gamma}{\partial\langle E\rangle}\langle I\rangle\\
    &\quad+
    \frac{\partial\log{Z}}{\partial\beta}\frac{\partial\beta}{\partial\langle E\rangle}
    +
    \frac{\partial\log{Z}}{\partial\gamma}\frac{\partial\gamma}{\partial\langle E\rangle}\\
    &=
    \beta
    \equiv
    \frac{1}{T},
  \end{split}
\end{equation}
where $T$ is called the temperature associated with the vortex configuration, and $\beta$ is the corresponding inverse temperature, or ``coldness'' according to \citet{garrod95}. From now on we will use the term
``temperature'' in this particular sense. By a process similar to \eqref{beq1T},
we can also obtain $\dfrac{\partial S}{\partial\langle I\rangle}=\gamma$, etc.~in case of more constraints.

\subsection{Three-dimensional vortex gas model}
A model of three-dimensional vortex gases is much more difficult and has been developed only in special cases. Giv\-en a vorticity field ${\boldsymbol\omega}({\bf x})$ in $\mathbb{R}^3$, the system
\begin{equation*}
  {\boldsymbol \omega}
  =
  \nabla\times{\bf V},
  \qquad
  \nabla\cdot{\bf V}
  =
  0
\end{equation*}
can be solved for the velocity field $\bf V(x)$ under the assumption that ${\boldsymbol\omega}({\bf x})$ decays sufficiently fast as $\|{\bf x}\|\to\infty$ (\citet{majda01}).
For the velocity field we get
\begin{equation*}
  {\bf V(x)}
  =
  -\frac{1}{4\pi}\int\frac{({\bf x-x'})\times{\boldsymbol\omega}({\bf x'})}{\|{\bf x-x'}\|^3}\,d\bf x',
\end{equation*}
and the kinetic energy can be written as
\begin{equation}
  \label{E}
  E
  =
  \frac{1}{8\pi}\iint\frac{{\boldsymbol\omega}({\bf x})\cdot{\boldsymbol\omega}({\bf x'})}{\|{\bf x-x'}\|}\,d{\bf x'}\,d{\bf x}.
\end{equation}

\subsubsection{Chorin's Monte Carlo Setup}
Vortex gases in three dimensions have been modeled in \citet{chorin-akao91} and \citet{chorin} assuming that vortices were made up of vertical and horizontal line segments connecting adjacent points on an integer lattice in the three-dimensional space, $\ZZ^3$. The vortices were supported on an oriented, self-avoiding random walk on the lattice. In the model, the expression \eqref{E} for the kinetic energy of one vortex becomes
\begin{equation*}
  E
  =
  \frac{1}{8\pi}\sum_i\sum_{j\ne i}\frac{{\boldsymbol\omega_i}\cdot{\boldsymbol\omega_j}}{\|i-j\|}
  +
  \frac{1}{8\pi}\sum_ iE_{ii},
\end{equation*}
where $i$ and $j$ are the three-dimensional coordinates of the locations of the centers of the vortex segments making up the vortex, $E_{ii}$ is the ``self-energy'' term of the $i$th segment that is a constant and is, therefore, neglected, $\|i-j\|$ is the distance between the $i$th and $j$th segment, and $\boldsymbol\omega_i$ is the vorticity of the $i$th segment. The velocity field was carefully defined using a cut-off function to prevent singularities and keep the terms in the sum defined. Similar ideas can be also used in the two-dimensional case.

To describe the state of such a vortex from a statistical viewpoint, \citet{chorin} defines the probability of a vortex with energy $E$ and inverse temperature $\beta$ via $P(E)=\dfrac{e^{-{\beta E}}}{Z}$, where $Z$ is a partition function. This probability is then maximized with respect to permissible configurations of the vortex for various fixed values of $\beta$ (positive, negative, and zero) using a Metropolis rejection algorithm. Note that from \eqref{beq1T}, as $\beta$ decreases from $+\infty$ to $0$, and then to $-\infty$, the temperature $T=1/\beta$ increases from $0$ through positive values to $+\infty=-\infty$ (denoted by $\infty$ in what follows), and then increases through negative values to $-0$.  To estimate the probability of various configurations Chorin considered ensembles of identical models and allowed them to evolve using the Metropolis algorithm.  Then Chorin estimated probabilities similar to the way one forecasts the weather using ensembles of weather models.  This does not take a large number of vortices.  The main information gained from the model is about qualitative behavior of the vortex, we are not suggesting that vortices are on a lattice.  However we think that this model provides useful qualitative information about the energy transfer from the outcome of the rather rapid dissipation that a vortex undergoes as it is stretched.

Chorin studied vortices of various lengths, and concluded that a stretched negative-temperature ($\beta<0$) vortex must fold assuming energy is conserved. In addition, as vortices with negative temperatures stretch, their temperature decreases to $T=\infty$ ($\beta=0$), and they become fractal. Note that as temperature crosses from negative to positive through $T=\infty$, the vortices tangle into three-dimensional configurations (see \citet{chorin-akao91,chorin}). It is also possible to model a fluid made up of several sparsely distributed vortices, but not without a significant increase in the computational complexity (\citet{marchioro94}). The experimental probabilities of the various configurations of the vortices can be found by performing Markov Chain Monte Carlo simulations (see \citet{chorin-akao91}).

\subsubsection{Entropy and Temperature}

A physical system is said to have negative temperature if when energy is added to the system the entropy of the system decreases, that is the system becomes less random.  We can give an example of such a system in fluid dynamics.  The supercritical vortex below a breakdown bubble in a vortex undergoing vortex breakdown is an example of such a system.  Two parameters have been used to study such vortices in a Ward chamber: the volumetric flow of the updraft and the swirl or angular momentum added to the flow.  It has been shown that when the ratio of the swirl to the volumetric flow of the updraft is in a certain range the flow configuration takes on certain structure.  This resembles a side-view of a champagne glass: the stem of the glass would correspond to the supercritical (what we will show to be the   negative temperature vortex) vortex and the breakdown bubble corresponds to the part of the glass above.

Theoretical studies by \citet{barcilon67,Burggraf77}, and \citet{fiedlerrotunno86} and experimental studies by \citet{ward72}, and later \citet{church77}, resulted in the identification of relationships between the radius of the supercritical vortex and the amount of swirl added to the flow, as well as relationships between the angular momentum added to the flow and the azimuth velocity and the vertical component of the velocity in the supercritical part of the vortex.  

Experimentally the following has been found: Holding the volumetric flow upward in the vortex chamber constant, as angular momentum is increased, the radius of the supercritical vortex decreases, as does the length of the supercritical vortex, while the vertical component of the velocity and the azimuthal component of the velocity increase.  This suggests that when the angular momentum is increased, the vortex is becoming less random as its volume decreases (entropy is decreasing) and the energy of the vortex is increasing as the velocity increases.  Hence the ratio of the change in entropy to the change in energy is negative and the temperature of the vortex is negative.   From the point of view of statistical mechanics such vortices transfer momentum to the larger scale flow when they dissipate.  This can happen in a couple of different ways. One way is by stretching (until the vortex kinks up and dissipates) and the other is by adding swirl until vortex breakdown occurs.  Video of tornados appears to show both possibilities occurring.  

Thus, negative temperature vortices form in vortex chambers when swirl and stretching are in a certain range.  As in the paper by \citet{fiedlerrotunno86} who argue that tornados can undergo the above conditions have similar behavior to vortices within vortex chambers, we believe that under similar conditions supercritical vortices in nature can also have negative temperatures.   The critical fact here is that the vortices must be in contact with the ground for the effect we are describing to occur i.e. for the negative temperature vortex to form.  The effects of the boundary layer are critical in that as the swirl increases the thickness of the boundary layer decreases, and the radius of the supercritical vortex (thought of as an extension of the boundary layer) also decreases as well as the length of the vortex. On the other hand, as the swirl increases the vertical velocity in the core and azimuthal velocity increase suggesting that the energy of the vortex increases.  We emphasize here that the swirl ratio must be in a certain range for the vortex to be supercritical, for a swirl ratio outside of this range the vortex behavior is different, and the temperature would not necessarily be negative.

We can think of the entropy of a collection of vortices (vortex gas) as a sum of its configurational entropy and a structural entropy corresponding to stretching, kinking, and collapse of each vortex. The structural entropy is described in the preceding paragraph.  Vortices moving into a region where there is a convergent updraft and an appropriate amount of low-level swirl will become more spatially organized (i.e., the configurational entropy decreases; see Figure \ref{Vortex sheet roll-up along the rear flank gust front}).

 We now argue that the vertical axis of vortices with negative temperature are straight and those with infinite temperature are fractal. Let us assume that $T<0$, so that also $\beta<0$. Since the probability of a configuration of a vortex with energy $E$ is $P(E)=e^{-\beta E}/Z$, the most likely configurations are those with a large positive energy $E$.

For a vortex with a small cross-section and nearly constant cross-sectional vorticity, the energy \eqref{E} will be maximized if the dot products, ${\boldsymbol\omega}({\bf x})\cdot{\boldsymbol\omega}({\bf x'})$, are as large as possible, hence the vortex should be straight. These comments can be generalized to a configuration with more than one vortex in which case they would be nearly parallel. Therefore, the highest-temperature (hottest) vortices (negative temperatures near $-0$) are the straightest. On the other hand, as $\beta=1/T$ approaches $-0$ (i.e., $T$ approaches $-\infty$), the probability distribution becomes uniform. In that case, all the orientations are equally likely, and the vortex configuration can be fractal. Hence, as negative-temperature vortices are stretched, they cool down and kink up as they dissipate.

Next we show, following \citet{chorin}, that if $T_1$ is the temperature of a vortex with mean energy $\langle E_1\rangle$, and $T_2$ is the temperature of a vortex with mean energy $\langle E_2\rangle$, and the vortex systems are combined, then, assuming conservation of energy, three cases are possible:
\begin{enumerate}
  \item If $T_2>T_1>0$, then $\dfrac{d\langle E_2\rangle}{dt}<0 $ and $\dfrac{d\langle E_1\rangle}{dt}>0$;
  \item If $T_1>0>T_2$, then $\dfrac{d\langle E_2\rangle}{dt}<0 $ and $\dfrac{d\langle E_1\rangle}{dt}>0$;
  \item If $0>T_2>T_1$, then $\dfrac{d\langle E_2\rangle}{dt}<0 $ and $\dfrac{d\langle E_1\rangle}{dt}>0$.
\end{enumerate}
Hence, for a vortex with a negative temperature, the closer the temperature is to $0$, the ``hotter'' the vortex is. If two vortices interact as a combined system, in which the hotter vortex will lose energy to the lower-temperature vortex, the system moves to an equilibrium state.

To show the three claims above, consider two vortices, one with energy $\langle E_1\rangle$ and a developing tornado with energy $\langle E_2\rangle$. We may regard this as two disjoint vortex systems, each separately in equilibrium. Assume the vortex with energy $\langle E_1\rangle$ moves into a developing tornado with energy $\langle E_2\rangle$ and that their probability densities are independent. Then the energy, $\langle E\rangle$, and the entropy, $\langle S\rangle$, of the combined vortex system are
\begin{equation*}
  \langle E \rangle
  =
  \langle E_1\rangle
  +
  \langle E_2\rangle
  \qquad\text{and}\qquad
  S
  =
  S_1+S_2.
\end{equation*}
As the combined vortex system adjusts to equilibrium, the time rate of change of the total entropy satisfies
\begin{equation*}
  \frac{dS}{dt}
  =
  \frac{dS_1}{dt}
  +
  \frac{dS_2}{dt}
  =
  \frac{dS_1}{d\langle E_1\rangle}\frac{d\langle E_1\rangle}{dt}
  +
  \frac{dS_2}{d\langle E_2\rangle}\frac{d\langle E_2\rangle}{dt}
  >
  0.
\end{equation*}
Conservation of energy implies that
\begin{equation*}
  \frac{d\langle E\rangle}{dt}
  =
  \frac{d\langle E_1\rangle}{dt}
  +
  \frac{d\langle E_2\rangle}{dt}
  =
  0.
\end{equation*}
Hence,
\begin{equation*}
  \frac{dS}{dt}
  =
  \left[
    \frac{dS_1}{d\langle E_1\rangle}
    -
    \frac{dS_2}{d\langle E_2\rangle}
  \right]
  \frac{d\langle E_1\rangle}{dt}
  =
  \left[
    \frac{1}{T_1}
    -
    \frac{1}{T_2}
  \right]
  \frac{d\langle E_1\rangle}{dt}
  >
  0,
\end{equation*}
and the three cases i--iii above follow.

Therefore, negative temperatures are ``warmer'' than positive temperatures, and negative temperatures that are closer to zero are warmer than negative temperatures that are farther from zero. Hence, as hot, negative-temperature vortices move into the developing tornado, they stretch and cool down, and the ambient vortex heats up. As this process repeats itself many times, the tornado vortex achieves a quasi-equilibrium with the environment.

Generalizing Chorin's argument used in the two-di\-men\-sion\-al model in \citet{chorin}, we consider a three-dimension\-al system in equilibrium, partition the system into boxes, and assume that the vortices are sparsely located, one per box, and nearly parallel and vertical. Let $E_{tot}=\frac12\sum_im_iU_i^2+\sum_i\langle E_i\rangle$, where $m_i$ is the mass of the $i$th box in the partition, and $U_i$ and $\langle E_i\rangle$ are the velocity and potential energy of the vortex in the $i$th box, respectively. If $T_i<0$, then $dS_i/d\langle E_i\rangle<0$, and hence, as the entropy increases, the energy is transferred from $\langle E_i\rangle$ to $U_i$. This explains the increase in the vorticity of the tornado as hot (negative-temperature) vortices enter the tornado. The vortices entering the tornado are stretched and begin to cool down and they kink up (In this situation the stretching dominates the swirl and the swirl ratio decreases.). As this happens, the entropy increases ($\Delta S>0$) and energy of the vortex decreases ($\Delta\langle E\rangle<0$) as it is transferred to the larger-scale flow, increasing the kinetic energy. This is the opposite of the situation considered above, where $\Delta S<0$ and $\Delta\langle E\rangle>0$, but $\beta=dS/d\langle E\rangle<0$.

We can illustrate a train of vortices in a vortex sheet entering the tornado and transferring the energy to the large scale by replacing the stretching term in the vertical vorticity evolution equation \eqref{vert-vort} by a train of delta functions in the form
 $\sum_{i}  \Gamma_i \delta_i =  \sum_{i} \Gamma_i \delta (t - i T)$ which entails

%

\begin{equation*}
  \frac{\partial\zeta}{\partial t}
  =
  -{{\bf V}\cdot\nabla\zeta}
  +
  \sum_{i} \Gamma_i \delta (t - i T)
  +
  {\boldsymbol\omega}_H\cdot{\nabla_Hw},
\end{equation*}
where $T$ is the time interval between successive vortices entering the tornado.
The response of the vertical vorticity to this forcing is a step function with jumps of the size $\Gamma_i$.  Hence the rate at which the vortices enter the ambient tornado vortex and their strength, and the rate at which the vorticity is removed from the ambient tornado vortex determine the eventual strength of the tornado.

\subsection{Energy Spectrum and Power Laws of Cai and Wurman}
\citet{chorin} gives two possible power laws of dissipation of energy with scale: $\langle E\rangle(k)\sim k^{-5/3}$ and $\langle E\rangle(k)\sim k^{-2}$, where $k$ is the wave number. The former law is derived by a scaling argument due to Kolmogorov and is also supported by Chorin's filament model using results from a Monte Carlo simulation. The latter law is derived as an alternative and is based on a possible form of the energy cascade. Chorin also claims that the latter one is a ``better candidate for the mean field result.''

Chorin's filament model
can be applied to analyze a vortex tube in a sparse, homogeneous suspension of tubes. Consider a narrow and straight enough vortex tube $\mathcal{T}$ that can be uniquely described by its center line, $\mathcal{C}$, parametrized by $s$, and cross-sections through $\mathcal{C}(s)$, denoted by $\mathcal{S}(s)$, orthogonal to the center line and such that $\mathcal{S}(s_1)$ and $\mathcal{S}(s_2)$ do not intersect for $s_1\ne s_2$.

Given a point ${\bf x}$ in the vortex tube and $r>0$, we define the ball $B_r({\bf x})=\{{\bf x'}\colon|{\bf x'}-{\bf x}|<r\}$. We take $r$ small enough so that $B_r({\bf x})$ contains no points that belong to other vortex tubes in the suspension. We denote by $\Sigma(s)$ the part of the cross-section $\mathcal{S}(s)$ inside $B_r({\bf x})$ and by $\mathcal{C}_r$ the part of the center line of the vortex tube for which $\Sigma(s)$ is non-empty, i.e.,
\begin{equation*}
  \Sigma(s)
  =
  \mathcal{S}(s)\cap B_r({\bf x})
  \quad\text{and}\quad
  \mathcal{C}_r
  =
  \left\{\mathcal{C}(s):\ \Sigma(s)\ne\emptyset\right\}.
\end{equation*}

In order to compute the energy spectrum, Chorin defines, for $r>0$, the vorticity correlation integral
\begin{equation*}
  S_r
  =
  \left\langle
    \int_{\mathcal{T}\cap B_r({\bf x})}\boldsymbol\omega({\bf x})\cdot\boldsymbol\omega({\bf x'})\,d\mathcal{H}_\mathcal{T}
  \right\rangle,
\end{equation*}
where $d\mathcal{H}_\mathcal{T}$ denotes the appropriate Hausdorff measure,
related to the set capacity on $\mathcal{T}$,
and the average is taken over the ensemble of all possible configurations. Then, using disintegration of measure (see \citet{lschwartz76}), we have
\begin{equation*}
  S_r
  =
  \left\langle
    \boldsymbol\omega({\bf x })\cdot\int_{\mathcal{C}_r}\,ds\int_{\Sigma(s)}\boldsymbol\omega({\bf x'})\,d\mathcal{H}_\Sigma
  \right\rangle.
\end{equation*}
If vorticity is roughly uniform throughout the cross-section $\Sigma(s)$ so that $\boldsymbol\omega({\bf x'})\approx\boldsymbol\omega(s)$, and if $|\Sigma|(s)=\mathcal{H}_\Sigma(\Sigma(s))$ denotes the Hausdorff measure of $\Sigma(s)$, then
\begin{equation*}
  S_r
  \approx
  \left\langle
    \boldsymbol\omega({\bf x})\cdot\int_{\mathcal{C}_r}|\Sigma|(s)\boldsymbol\omega(s)\,ds
  \right\rangle.
\end{equation*}
Assuming that the Hausdorff dimension of $\Sigma(s)$ remains constant throughout $\mathcal{T}\cap B_r({\bf x})$ and denoting it by $D_{\Sigma}$, we obtain
\begin{equation*}
  S_r
  =
  \mathcal{O}(r^{D_{\Sigma}+1}).
\end{equation*}
To obtain the vorticity spectrum, $Z(k)$, we integrate the Fourier transform of $S_r$ over a sphere of radius $k=|{\bf k}|$ (\citet{chorin-akao91,chorin}). This gives $Z(k)=\mathcal{O}(k^{-D_{\Sigma}+1})$, and, consequently, the energy spectrum satisfies
\begin{equation}\label{e-k-ksq}
\langle E\rangle(k)=Z(k)/k^2=\mathcal{O}(k^{-D_{\Sigma}-1}).
\end{equation}

Let $D_c$ be the dimension of the center line $\mathcal{C}$ of the vortex, and $D$ the dimension of the support of the vorticity in the vortex filament $\mathcal{T}$. We argued above that for vortices with negative temperature the center line of the vortex has Hausdorff dimension one, i.e., $D_c=1$, and therefore (see \citet{pesin09})
\begin{equation*}
  D
  =
  D_c
  +
  D_\Sigma.
\end{equation*}

Supported by the results of \citet{fujita81,fiedlerrotunno86,fiedler94,lewellensxia00,wurman02,xialewellens03,lewellens07,lewellens07a},
we suggest that tornadoes have fractal cross-sections and subvortices moving into the larger tornadic flow have negative temperature. Graphical evidence is provided in Figures \ref{Vortex sheet roll-up along the rear flank gust front}, \ref{fig:suctionspots}, and \ref{fig:tornado_fractal}. This suggestion is further supported by the existence of subvortices within subvortices (see Figure \ref{fig:hier-vort}) as shown in recent videos (\citet{wadena10,elrenovideo13}).

Assume now that the cross-section of the tornado is fractal for $T<0$ with cross-sectional dimension $D_\Sigma$. For $ 1<D_\Sigma\le2 $ the energy (see \eqref{e-k-ksq}) satisfies $\langle E\rangle(k)=\mathcal{O}(k^{-\gamma})$ with $2<\gamma\le3$. It follows that for large scales (small $k$), an increase in $\gamma$ in the range from $2$ to $3$ corresponds to an increase of the energy $\langle E\rangle(k)$. This is consistent with the idea that vortices from a vortex sheet feeding a larger tornado cause an increase in the dimension of the tornado's cross-sectional area and an increase in its energy. Thus, an increase in the dimension of the cross-sectional area may be associated with tornadogenesis or strengthening of an existing tornado. Note that this is analogous to Cai's power law,
in which a decrease in the (negative) exponent in the power law \eqref{eq:cai_power_law} leads to tornadogenesis or a stronger tornado.

To gain further insight into the processes that might contribute to tornadogenesis, we consider the effects of helicity. First we consider its effects in the simplified case of homogeneous isotropic turbulence (see \citet{lilly86a,lilly86b}). Writing the Navier--Stokes equation in energy form and taking the Fourier transform, \citet{chorin} obtains
\begin{equation*}
  \partial_t\langle E\rangle(k)
  +
  2k^2R^{-1}\langle E\rangle(k)
  =
  Q(k),
\end{equation*}
where $Q(k)$ comes from the nonlinear term in the Navier--Stokes equation. This term represents the transfer of energy between wave numbers and has been studied extensively for the case of homogeneous turbulence (see, e.g., \citet{waleffe92}). Certain terms have been singled out and studied in relation to inverse energy cascades. These interactions involve three wave numbers. It was found that the net effect of the so-called nonlocal interactions is to transfer energy from intermediate scales to larger scales. These interactions occur between modes with helicity of the same sign. 
\citet{andre77} and \citet{lilly86b}
discuss flows without helicity and dissipation of energy as shown in Figure \ref{lilly-13}, and flows with helicity and low dissipation of energy as shown in Figure \ref{lilly-15}. The images represent the results of two numerical experiments and show that isotropic turbulence with helicity inhibits the dissipation of energy at large scales.

Under the effects of strong rotation, the flow has the tendency to become anisotropic (\citet{pouquet10}). Studies have shown that the presence of helicity and low-energy dissipation are unlinked unless the helicity is continuously supplied and/or generated at the energy-containing scales; this is associated with inhomogeneity in the mean field (\citet{lilly86a,lilly86b,yokoi93}). Such an inhomogeneity would be supplied by surface friction and the rear flank and forward flank downdrafts and/or their gust fronts. The increase in the exponents for the power laws for the vorticity as tornadogenesis approaches is consistent with the helicity production of the flow at the energy-containing scales. Idealized cross-sections of vortices with high helicity exhibit self-similarity (compare Figure \ref{fig:hier-vort} and the image in \citet[p.~168]{arnold99}). Thus the tornado develops at a ``focus'' scale at which much energy and helicity is transferred among other scales. This is consistent with the mean power law and helicity contribution to the flow (\citet{pouquet10}).


\section{Conclusions and Future Work}
\label{sec:conclusions}
The two-dimensional point vortex theory presented in Section \ref{sec:vortex_gases} can be applied to pairs of cyclonically rotating vortices in the half-plane. The paths of the pairs of interacting vortices (see Figure \ref{fig:f-pics}) form the same type of pattern as the tracks of overlapping suction vortices moving through fields (see Figure \ref{fig:suctionspots}) as observed by \citet{fujita81} and others
(\citet{grazulis97}) from the air. As noted earlier (cf.~Section \ref{sec:vortex_gases}), under the influence of strong rotation, turbulent flow becomes anisotropic with the flow tending toward, but never fully becoming, a two-dimensional flow (\citet{pouquet10}). This suggests that the two-dimensional flow may be an attractor in this situation. We suggest that the increase in the power in Wurman's and Cai's power laws is an indication of the anisotropy of the flow. The paths in Figure \ref{fig:suctionspots} can be modeled using the two-dimensional vortex gas theory with translating and interacting point vortices (Figure \ref{fig:f-pics}). In the first pair of paths in Figure \ref{fig:f-pics}, one can identify the two overlapping paths of cyclonically rotating vortices in the left half-plane and two overlapping paths of counter-rotating, anticyclonic, mirror vortices in the right half-plane as the other ends of the arching vortex lines. The image on the right in Figure \ref{fig:f-pics} shows the paths left by two translating interacting vortices. Both the first pair of paths and the last path in Figure \ref{fig:f-pics} look similar to paths in Figure \ref{fig:suctionspots}. The tracks left by the suction vortices (Figure \ref{fig:suctionspots}) suggest the vortices either originate in the larger tornado vortex or move into the tornado vortex, after which they intensify due to stretching, make a partial revolution, and then dissipate. In either case they appear to occur as a result of the roll-up of a vortex sheet as observed in numerical simulations (see, e.g., \citet{rotunno84} and Figure \ref{vortex_image_time_series1.pdf}). Some of these intense vortex paths are from one to two yards in diameter (Figure \ref{fig:suctionspots}), and some are as narrow as $30$ cm in diameter (\citet{fujita81}).
These vortices, which may be extremely intense, have been observed to pull cornstalks out of clay soil by their roots. We suggest that these vortices have negative temperature and transfer energy to the larger vortex as they dissipate. This could manifest itself as a vortex breakdown. Indeed, video footage of subvortices in tornadoes suggests that they behave as negative-temperature vortices would. For example, in some instances the vortices' appearance is associated with stretching and with strong convergence. This may indicate that the vortex intensification is related to a decrease in entropy and an increase in energy. Numerical simulations of intense vortices in \citet{fiedlerrotunno86,fiedler94,lewellensxia00,xialewellens03,lewellens07,lewellens07a} show that the maximum wind speeds in intense narrow vortices undergoing vortex breakdown may exceed the speed of sound in the vertical direction.

Mobile Doppler radar observations of the vortices apparently originating in the vortex core suggest that they make a partial revolution about the ambient tornado vortex and then dissipate. This might indicate a Hopf bifurcation of the horizontal component of the flow field, creating a two-cell flow structure: downdraft core and a surrounding updraft. \citet{wurman02} has found evidence supporting both the creation of vortices inside the tornado and outside the tornado resulting in flows potentially enhancing the tornado's strength. He has also noted that these secondary vortices have a different velocity and shear profile than the parent tornadoes. The tornadoes appear to have a two-cell structure and a modified Rankine combined profile, with mean azimuthal velocity, $v$, depending linearly on radius inside the tornado core, $v=Cr$, and a power-law outside the core, $v=Cr^{-b}$, where $0.5\le b\leq0.6$. In an extreme case, \citet{wurman02} found a power law $v=Cr^{-1}$ on one side of an intense tornado. Such a power law would be consistent with no vorticity outside the tornado core on that side. On the other side of the tornado,
the power law was found to be $v=Cr^{-b}$ with $0.5\le b\leq0.6$, which is consistent with vorticity being advected into the tornado, possibly along a vortex sheet spiraling into it. The secondary vortices appear to be single-celled with extreme values of shear and extreme transient updrafts. This is consistent with these vortices having negative temperatures in the vortex gas sense (\citet{chorin}).

\begin{figure*}
  \begin{center}
    \includegraphics[width=0.13\textwidth]{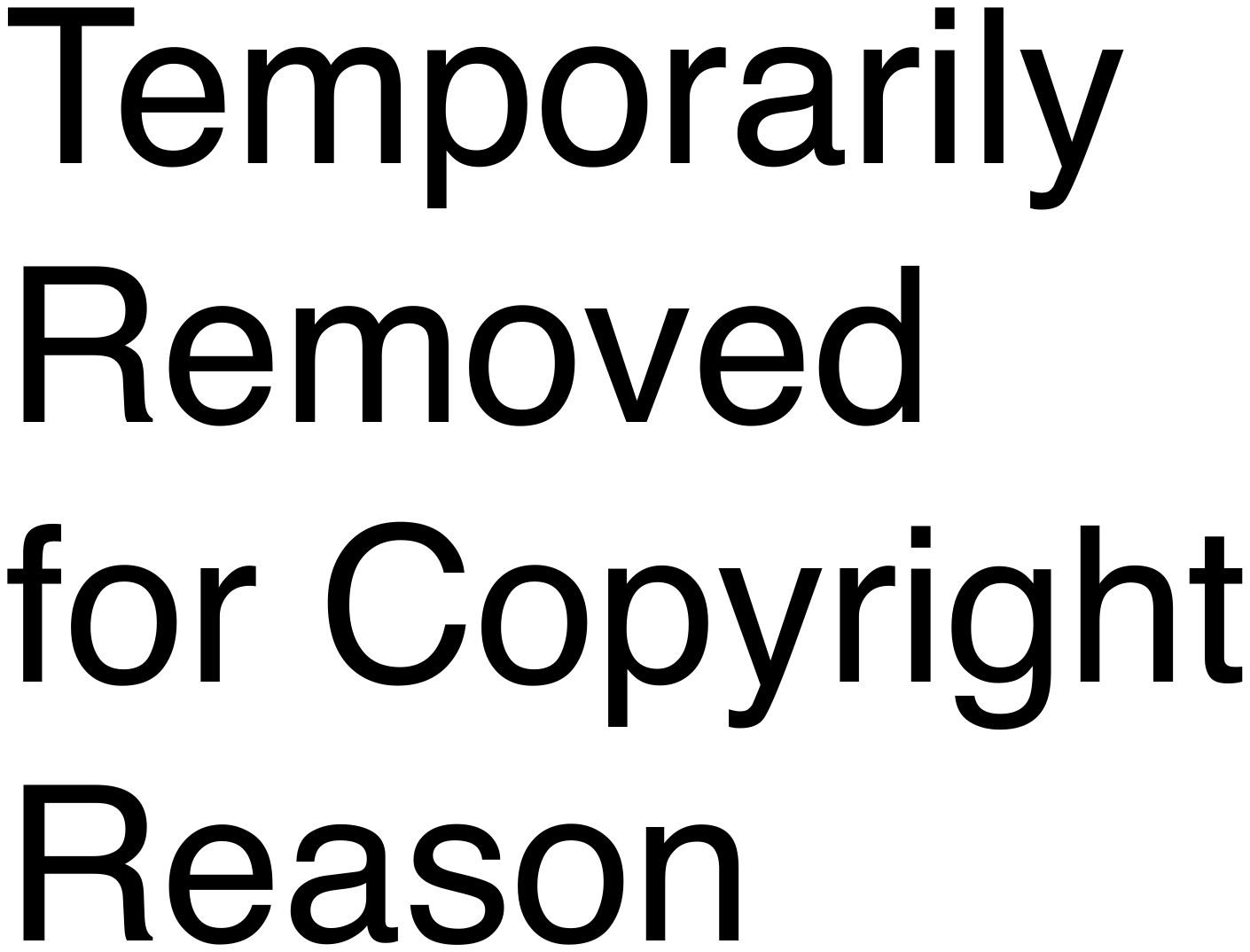}
  \end{center}
  \caption{Vortex sheet roll-up; \copyright~Gene Moore.}
  \label{Vortex sheet roll-up along the rear flank gust front}
\end{figure*}

\begin{figure*}
  \begin{center}
    \includegraphics[width=0.97\textwidth]{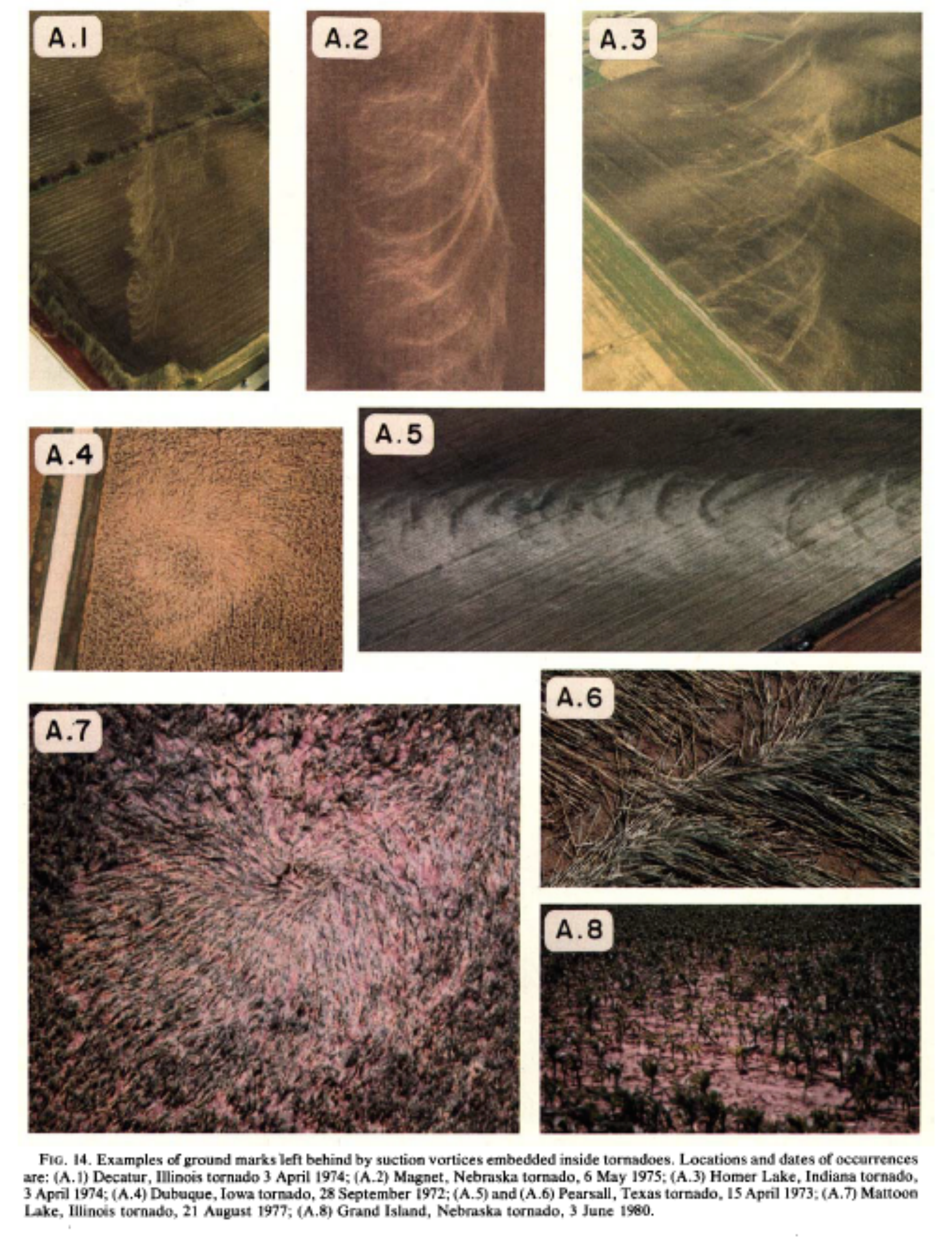}
  \end{center}
  \caption{Tracks left in corn fields showing vortices spiraling into tornadoes and then dissipating; \copyright~AMS, \citet{fujita81}.}
  \label{fig:suctionspots}
\end{figure*}

\begin{figure*}
  \begin{center}
    \includegraphics[width=0.89\textwidth]{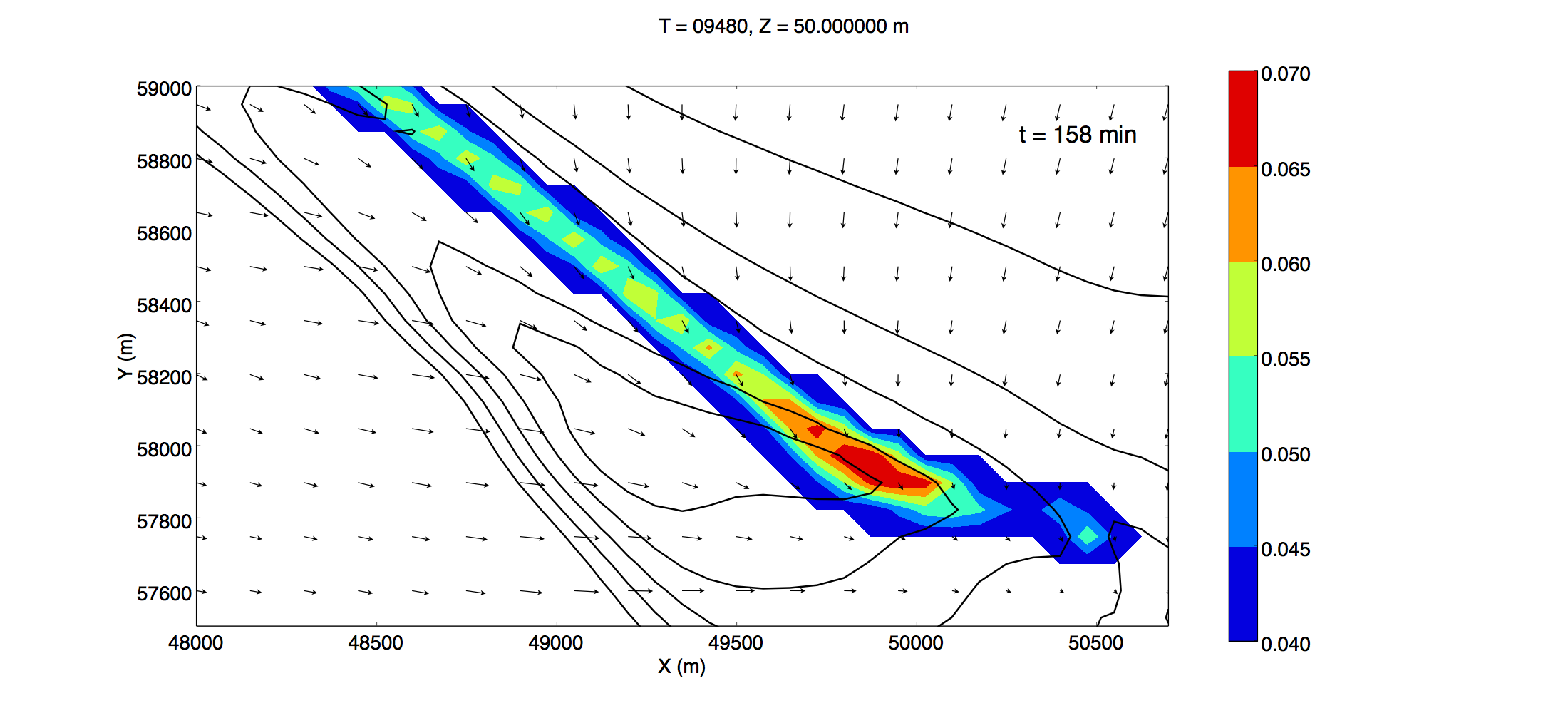}
    \includegraphics[width=0.89\textwidth]{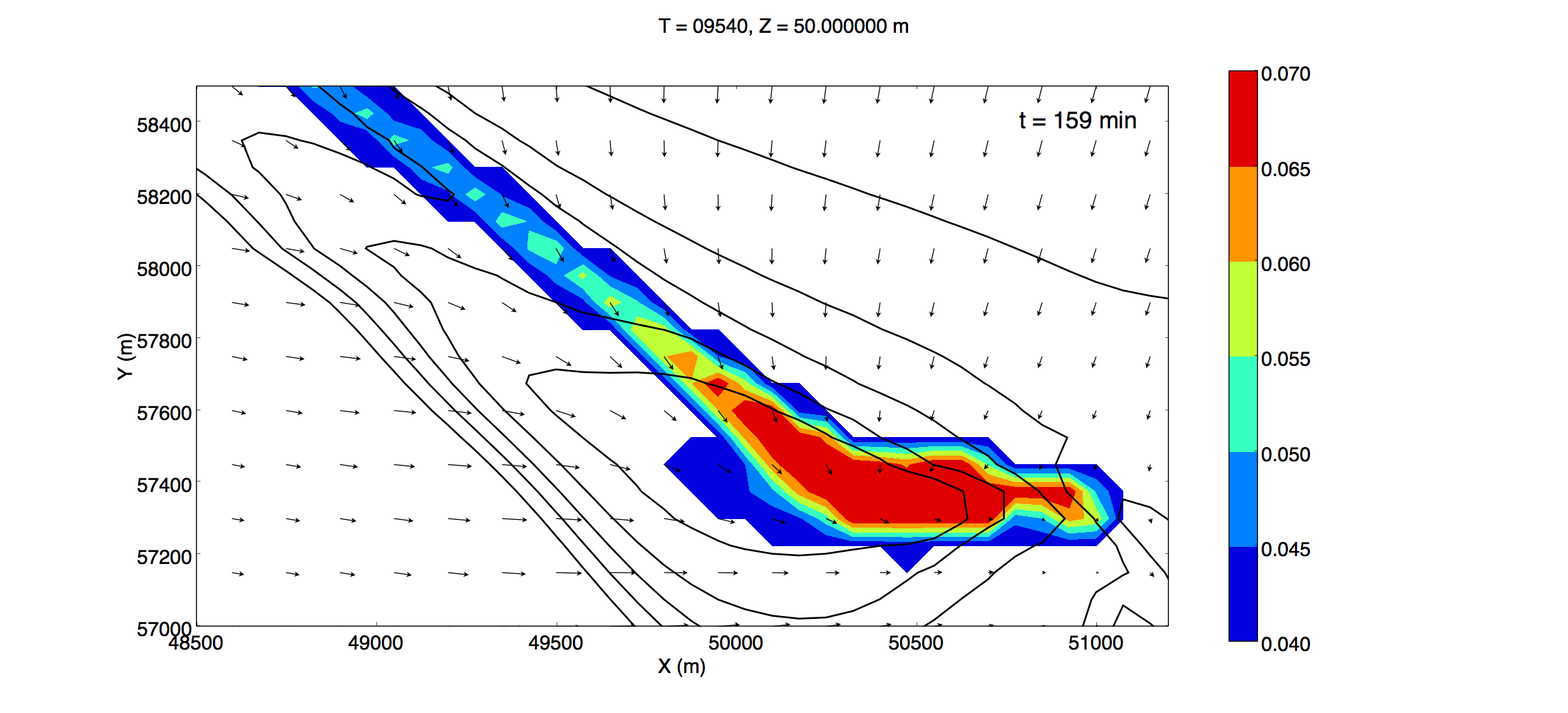}
    \includegraphics[width=0.89\textwidth]{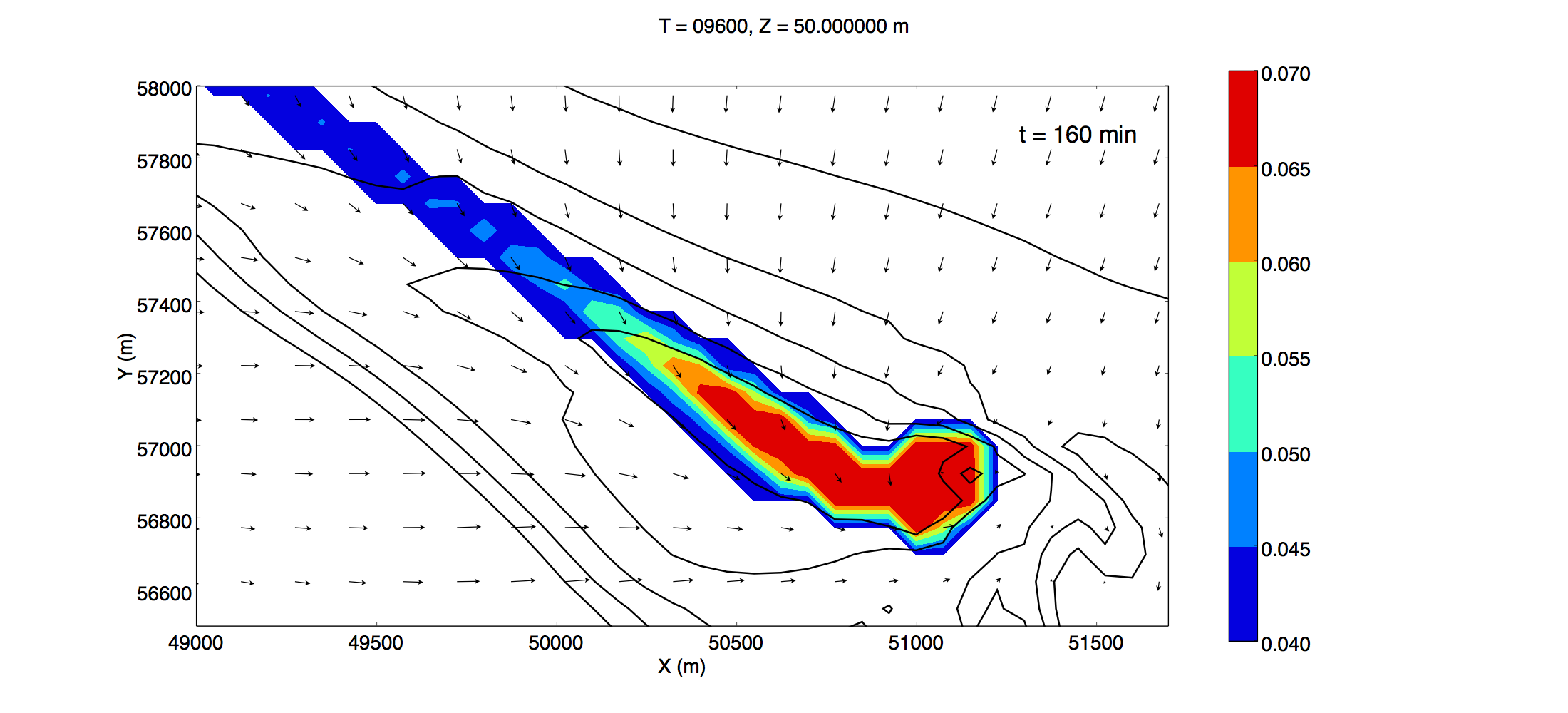}
  \end{center}
  \caption{Vortex sheet roll-up in CM1 simulation at $z=50$ meters and $t=158$, $159$, and $160$ minutes. Vertical vorticity is shaded, vertical velocity is contoured in $1$ m s$^{-1}$ intervals beginning at $1$ m s$^{-1}$, and horizontal wind vectors are plotted as arrows.} 
  \label{vortex_image_time_series1.pdf}
\end{figure*}

\begin{figure*}
  \begin{center}
    \includegraphics[width=0.66\textwidth]{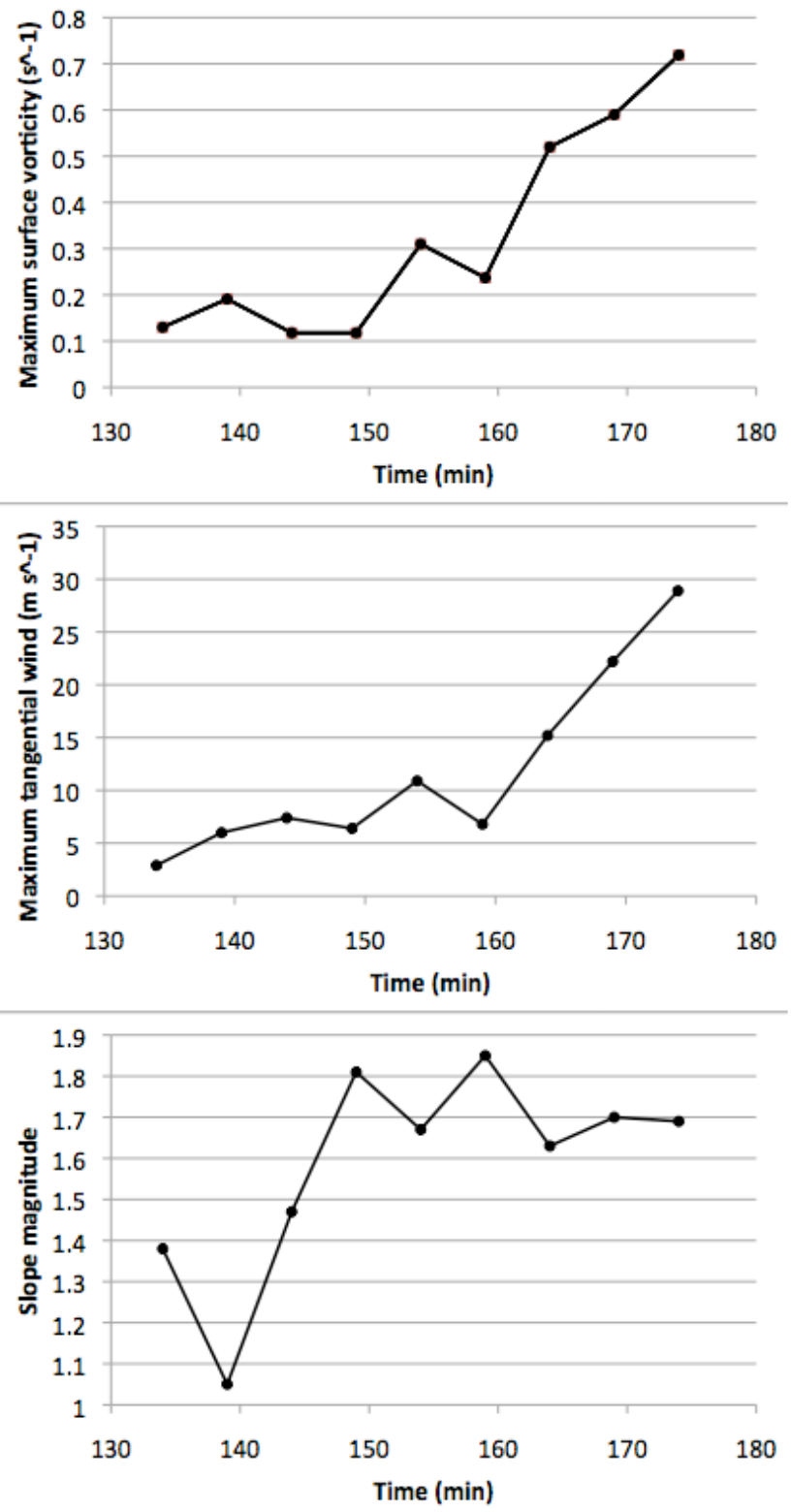}
  \end{center}
  \caption{Time series of (a) maximum surface vorticity (s$^{-1}$), (b) $V_T$ (m s$^{-1}$), and (c) vorticity line slope magnitude.}
  \label{timeseries-ok-nssl}
\end{figure*}

\begin{figure*}
  \begin{center}
    \includegraphics[width=0.85\textwidth]{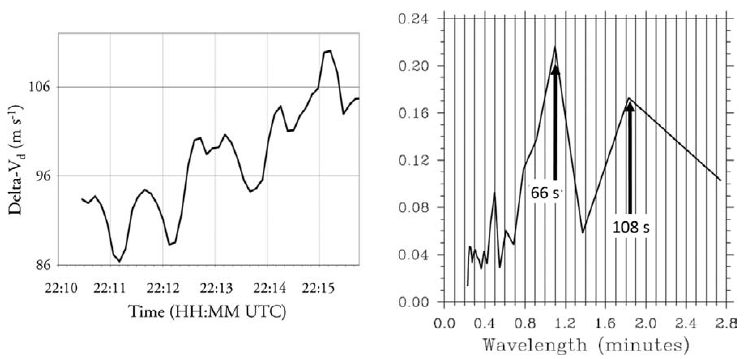}
  \end{center}
  \caption{Time series of maximum gate-to-gate shear, $\Delta V$, consistent with energy being pumped into the tornado in discrete pulses, possibly from roll-up vortices within a vortex sheet (left); FFT of $\Delta V$ with peak energy at $66$ seconds and $108$ seconds (right), suggestive of rotating asymmetry in the vortex; \copyright~AMS, \citet{wurman2013}.}
  \label{wurman-kosiba-fig7}
\end{figure*}

\begin{figure*}
  \begin{center}
    \includegraphics[width=0.73\textwidth]{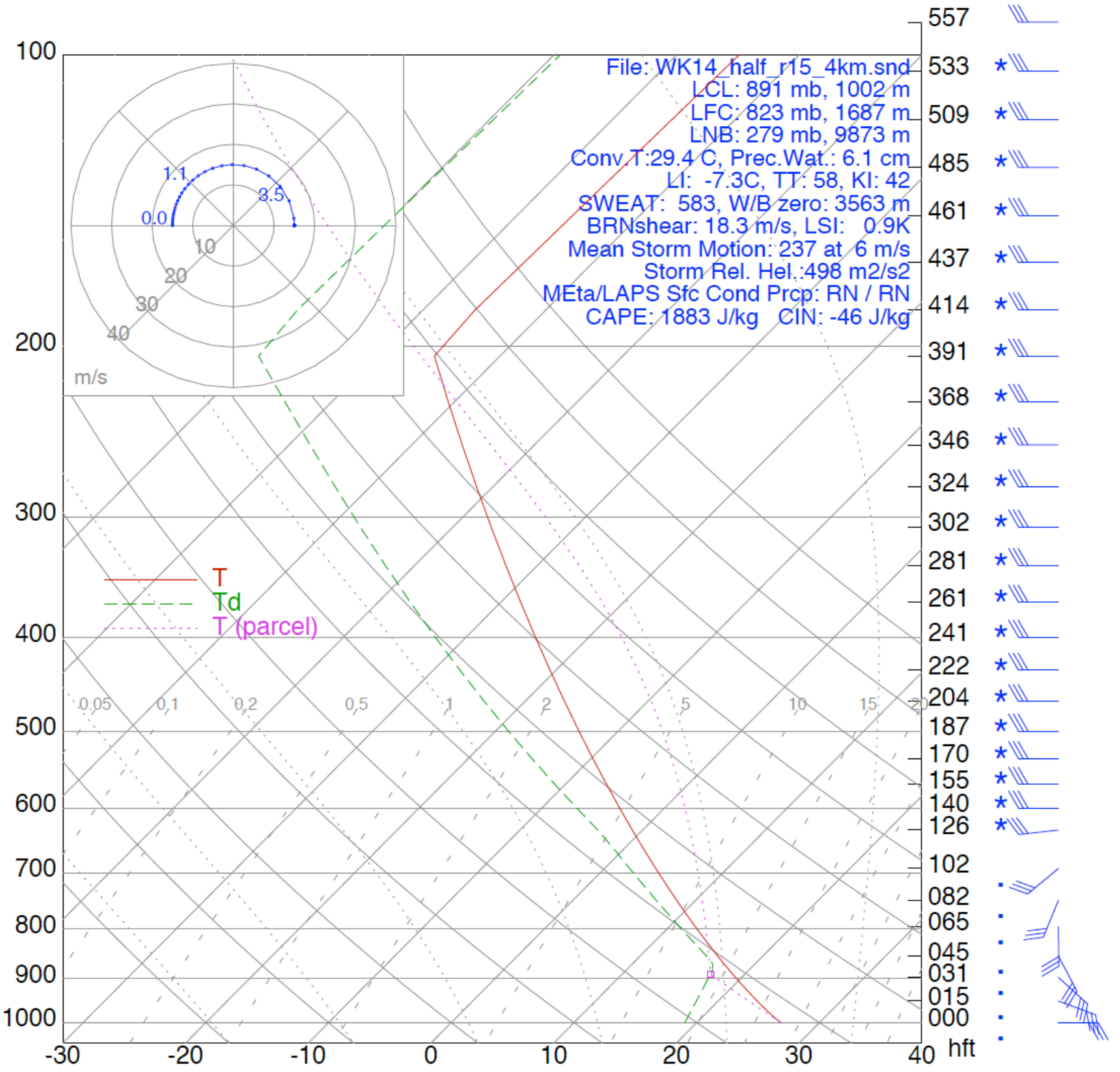}
  \end{center}
  \caption{Base state for CM1 simulation. Wind barbs are in knots. The hodograph (upper left corner) is in m s$^{-1}$, with marked heights in km.}
  \label{skewt-ok-nssl}
\end{figure*}

\begin{figure*}
  \begin{center}
    \includegraphics[width=0.71\textwidth]{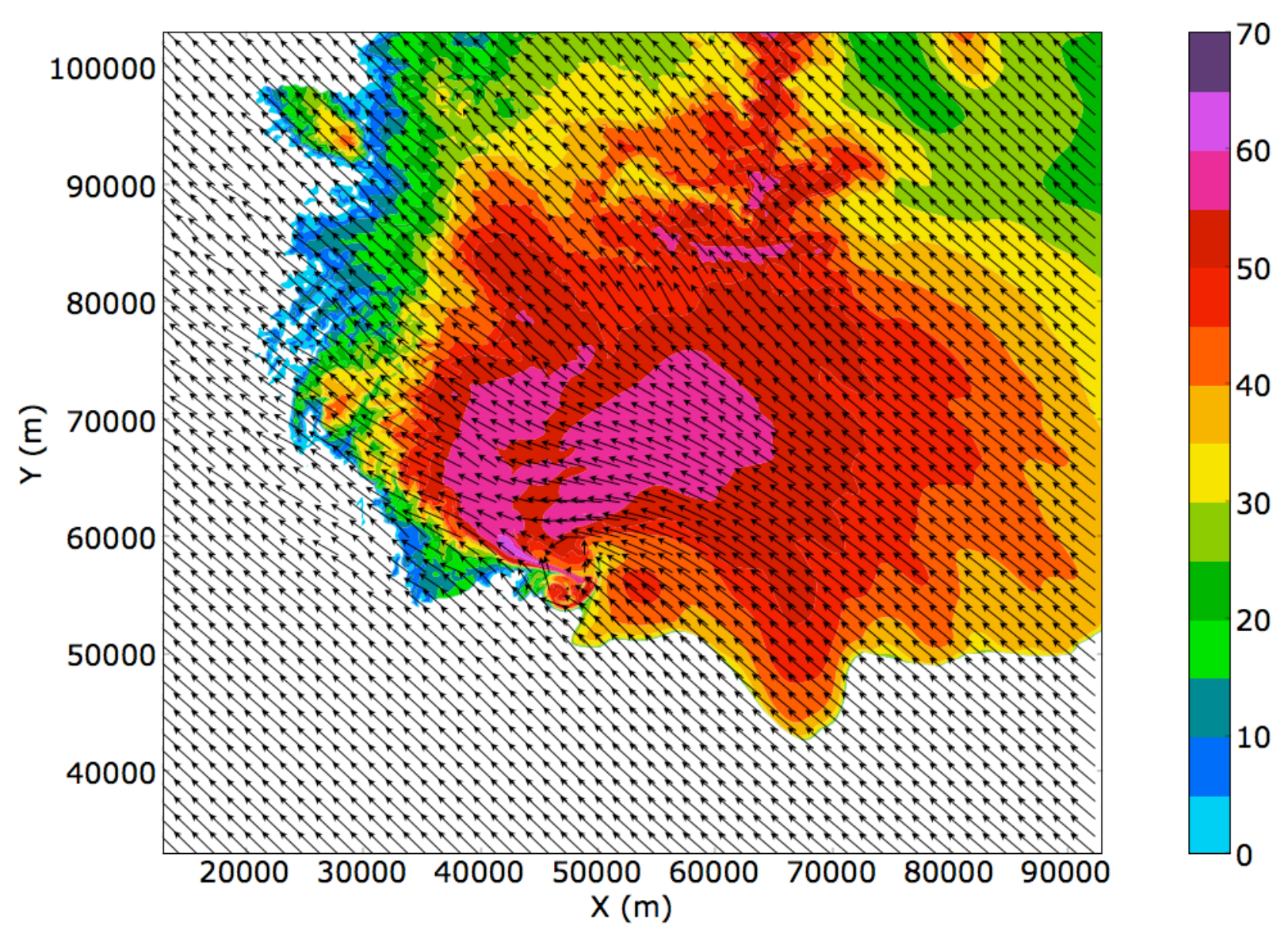}
  \end{center}
  \caption{Simulated reflectivity field (dBZ) and horizontal wind vectors (arrows) valid $\sim 0.5$ km AGL at $t=164$ min, near the time of tornadogenesis.}
  \label{simulated_fields}
\end{figure*}

\begin{figure*}
  \begin{center}
    \includegraphics[width=0.6\textwidth]{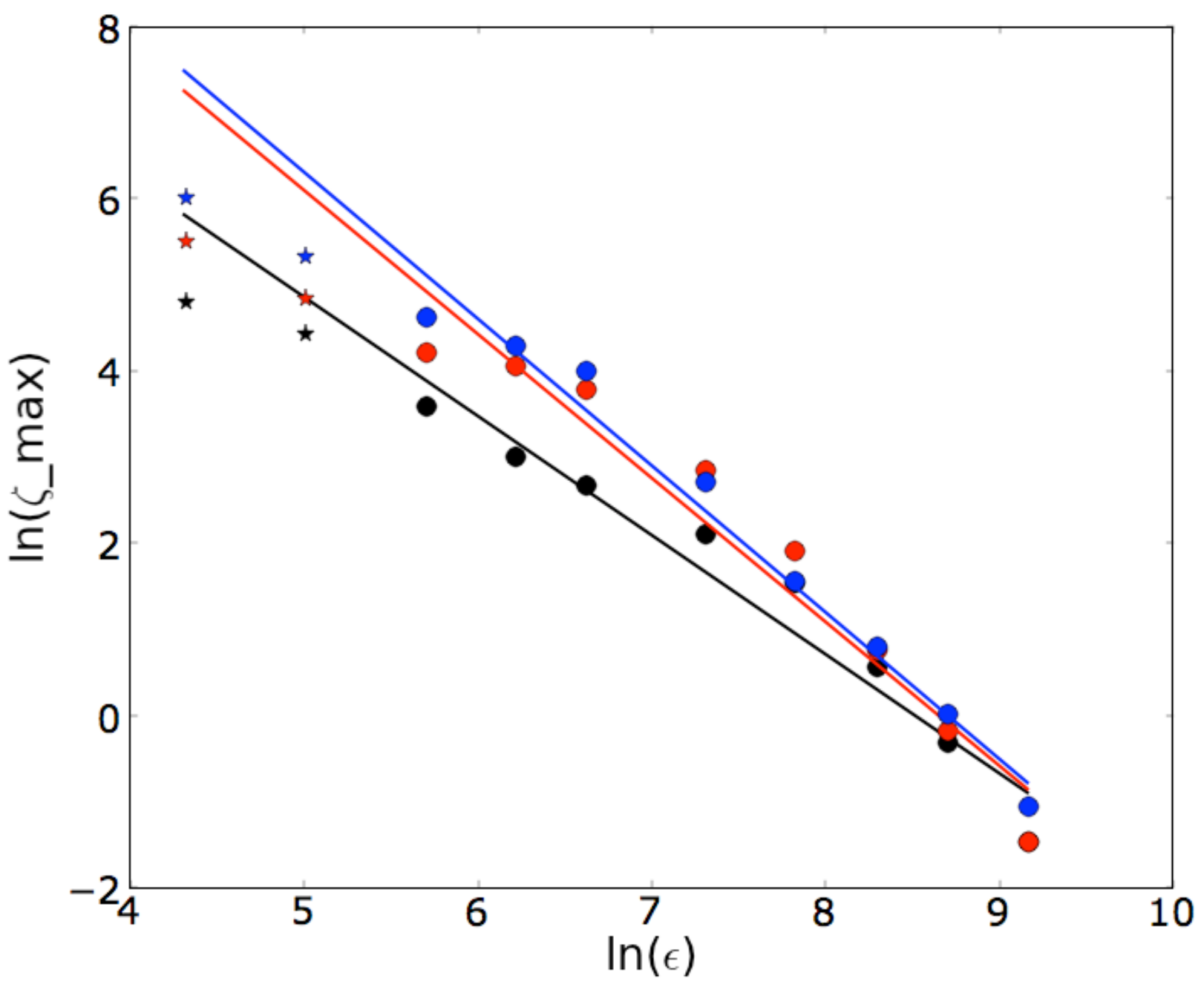}
  \end{center}
  \caption{Vorticity lines computed at $t=134$ min (black, least steep), prior to the development of a discernible surface vortex; at $t=154$ min (red, in-between), by which time a relatively weak surface vortex is present; and at $t=169$ min (blue, steepest), near the time of tornadogenesis. Points used to create the least-squares fit are denoted by dots, while points not used in the vorticity line computation ($\varepsilon<300$ m) are denoted by asterisks.}
  \label{vorticitylines-ok-nssl}
\end{figure*}

The results of Section \ref{sec:cm1}, the last part of Section \ref{sec:self-similarity}, and Section \ref{sec:vortex_gases} suggest that the increase in strength of a developing tornado occurs as a result of the increase in vorticity due to an inverse energy cascade from smaller scales. This process appears to be related to a negative viscosity phenomenon described by \citet{lilly69,lilly76} and \citet{bluestein13}. Photographic evidence supporting this is given in the photo by Gene Moore in Figure \ref{Vortex sheet roll-up along the rear flank gust front}. Numerical models have also supported this. A possible connection between tornadoes and nearly continuously (periodically) produced vortices (vortex lines (tubes)) is that the vortices stir or pump the tornado and increase the vorticity. The frequency with which vortex lines are produced, their strength, and the stretching of the vortices determine the eventual strength of the tornado. This can be seen from the point of view of the vortex gas theory (Section \ref{sec:vortex_gases}).

\begin{figure*}
  \begin{center}
    \includegraphics[height=0.3\textheight]{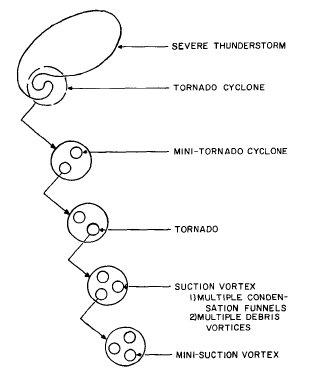}
    \hfill
    \includegraphics[width=0.55\textwidth]{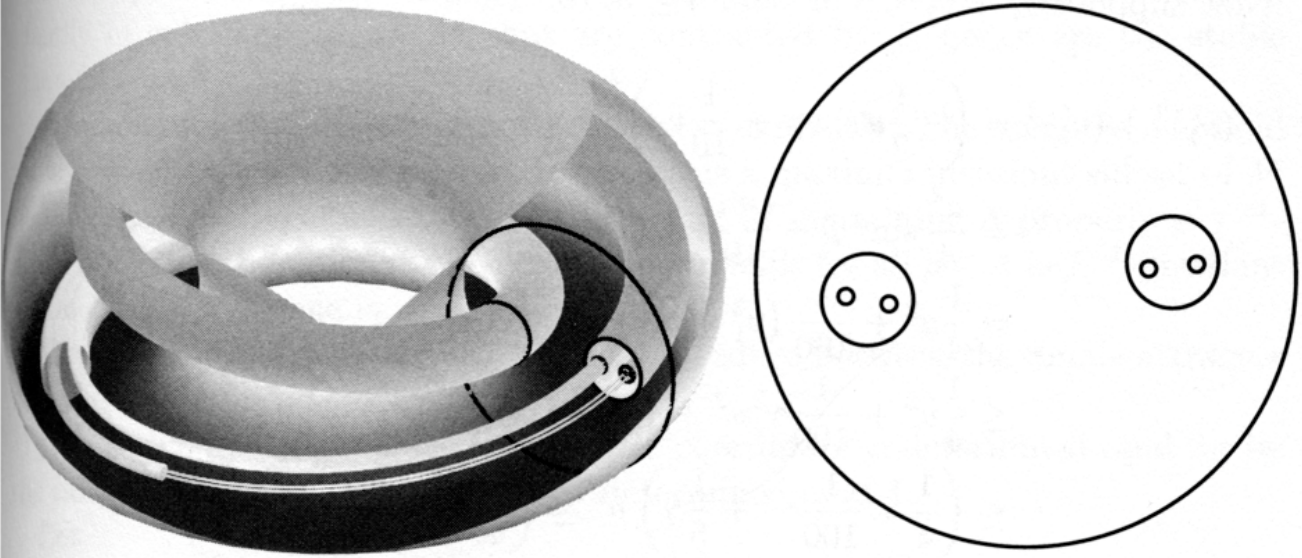}
  \end{center}
  \caption{Hierarchy of known vortex scales in tornadic supercells (left); \copyright~AMS, \citet{church77}. Smale--Williams attractor and its cross-section (right); \copyright~Cambridge University Press, \citet{katok-hasselblatt95}.}
  \label{fig:hier-vort}
\end{figure*}

\begin{figure*}
  \begin{center}
    \includegraphics[width=\textwidth]{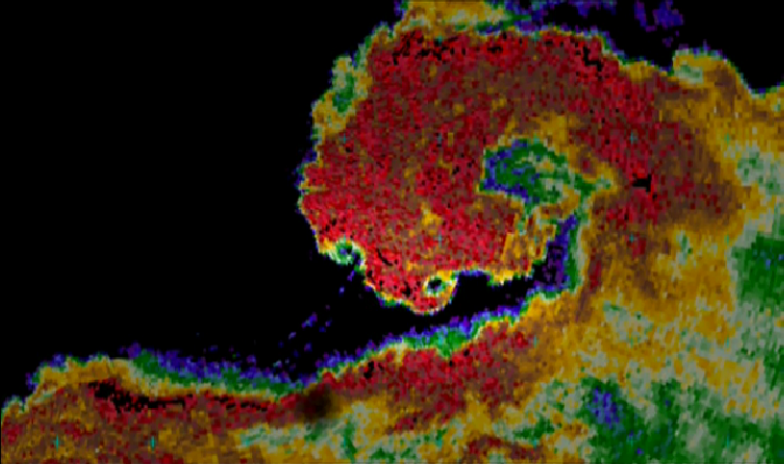}
  \end{center}
  \caption{A reflectivity image of a tornado showing self-similarity and a possible fractal structure; \copyright~Joshua Wurman, \citet{nova04}.}
  \label{fig:tornado_fractal}
\end{figure*}

\begin{figure*}
  \centering
  \includegraphics[height=0.35\textheight]{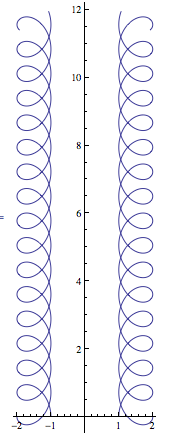}
  \hskip 1cm
  \includegraphics[height=0.35\textheight]{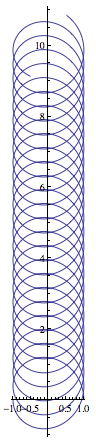}
  \caption{Graphs of interacting point vortices have tracks similar to suction vortex tracks.  The left image is similar to the images A.1 and A.2 of Figure \ref{fig:suctionspots}. The right image is similar to the image A.5 of Figure \ref{fig:suctionspots}.}
  \label{fig:f-pics}
\end{figure*}

\begin{figure*}
  \begin{center}
    \includegraphics[width=0.7\textwidth]{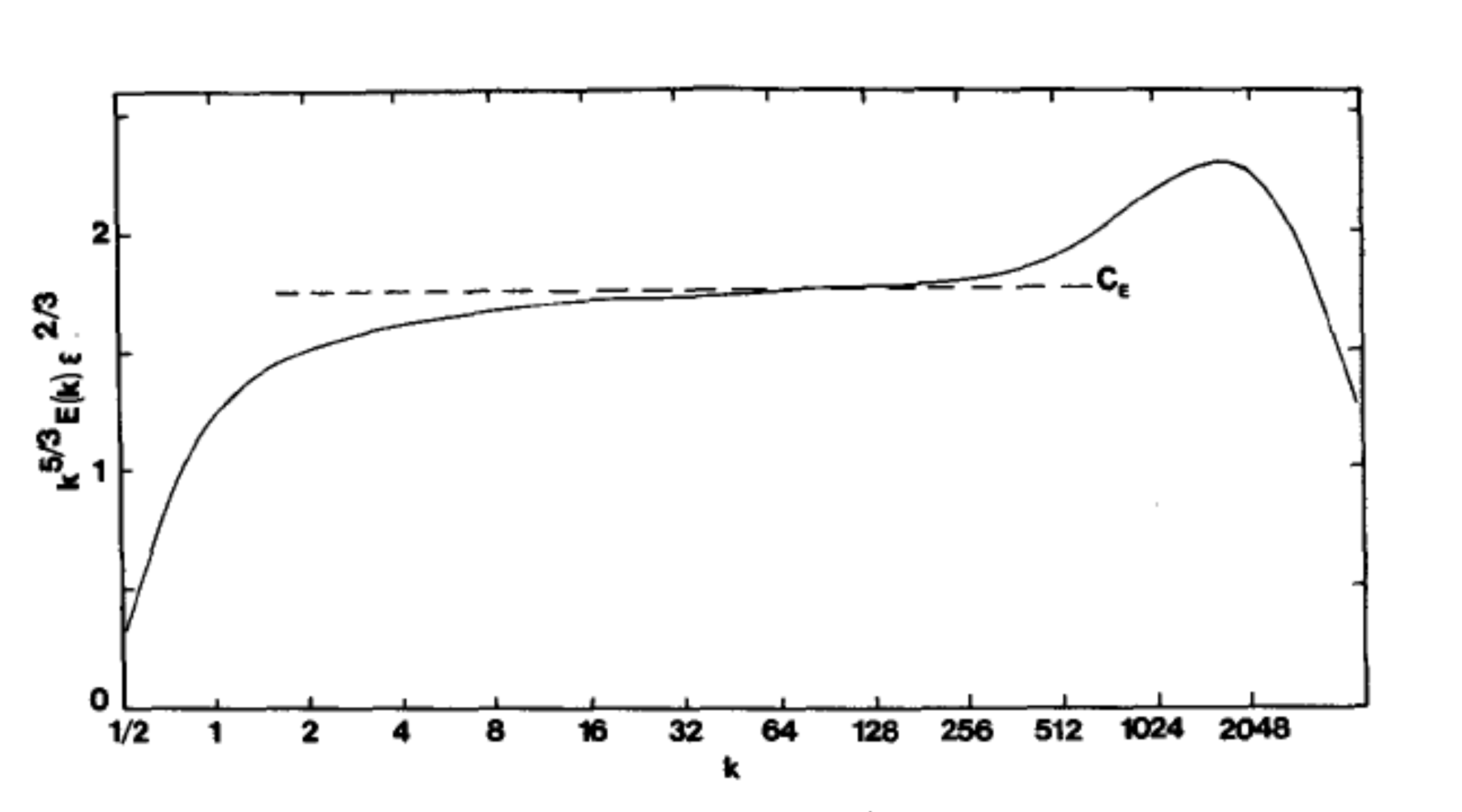}
  \end{center}
  \caption{Inertially weighted energy spectrum corresponding to time $t=8$ for a turbulent flow with low helicity showing high dissipation at large scales (small $k$); \copyright~AMS, \citet{andre77,lilly86b}.}
  \label{lilly-13}
\end{figure*}

\begin{figure*}
  \begin{center}
    \includegraphics[width=0.7\textwidth]{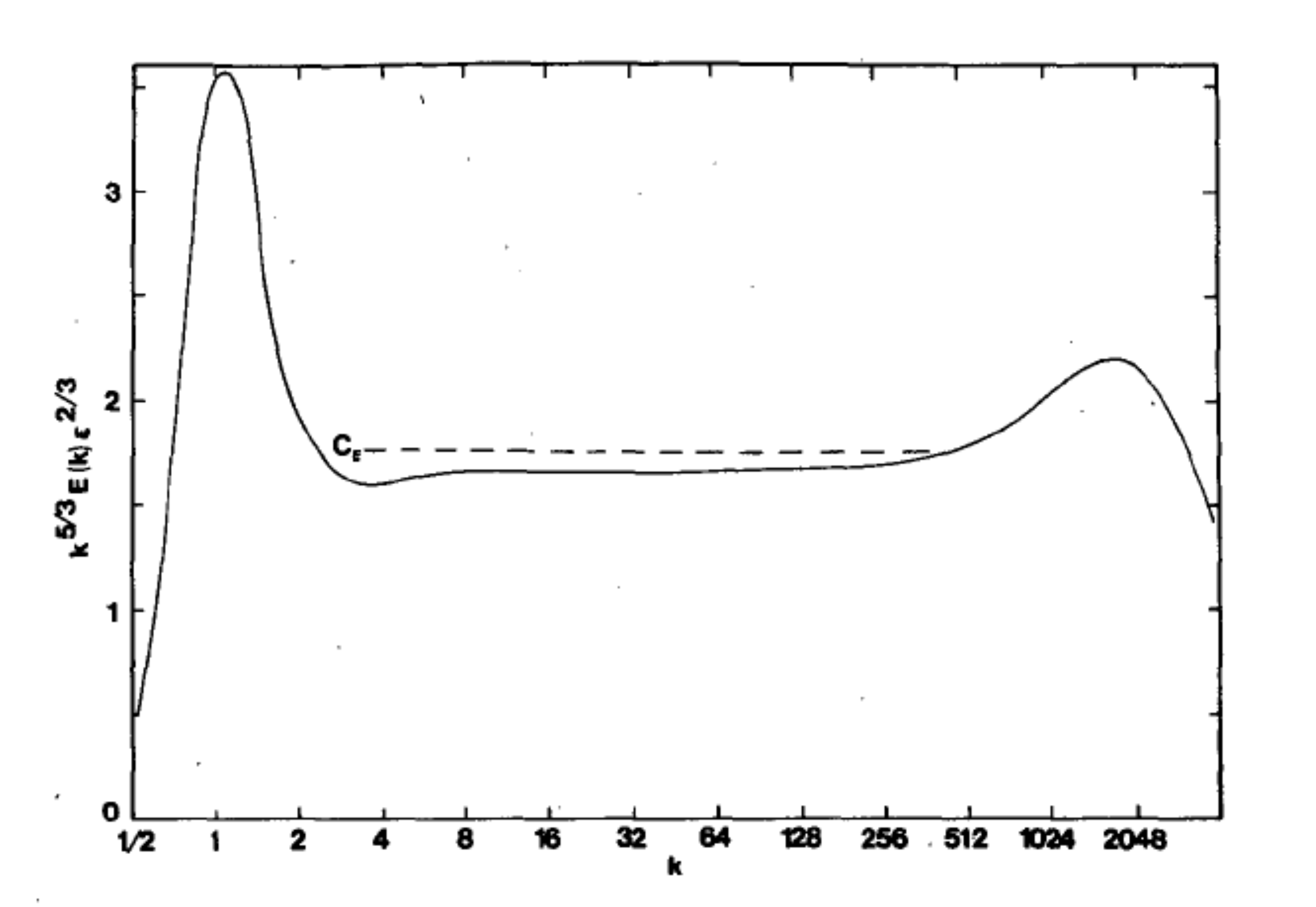}
  \end{center}
  \caption{Inertially weighted energy spectrum corresponding to time $t=12$ for a turbulent flow with high helicity showing low dissipation at large scales (small $k$); \copyright~AMS, \citet{andre77,lilly86b}.}
  \label{lilly-15}
\end{figure*}

Supercell thunderstorms have a quasi-periodic nature, cycling between destructive and rebuilding phases. They typically have a longer life span than generic storms. Curiously, most classic supercells have common features which make them distinctive from other storms of the same scale. These features include, for example, a wall cloud, a tail cloud, and a flanking line. This commonality of features suggests that atmospheric flows that demonstrate themselves as supercells fluctuate near ``attractor'' flows that share some common structure. In the radar reflectivity image of the hook echo region of a supercell thunderstorm shown in Figure \ref{fig:tornado_fractal}, the hooks on the boundary of the region represent successive vortices in a vortex sheet. Such vortices could provide periodic pulses of energy to the tornado. An analogy can thus be drawn with the work of \citet{kuznetsov11}, in which periodic pulses introduced into a dynamical system lead to a Smale--Williams attractor. As noted in Section \ref{sec:self-similarity}, the hierarchy of vortices shown in Figure \ref{fig:hier-vort} is reminiscent of the cross-sections of the Smale--Williams attractor.

\begin{figure*}
  \begin{center}
    \includegraphics[width=0.7\textwidth]{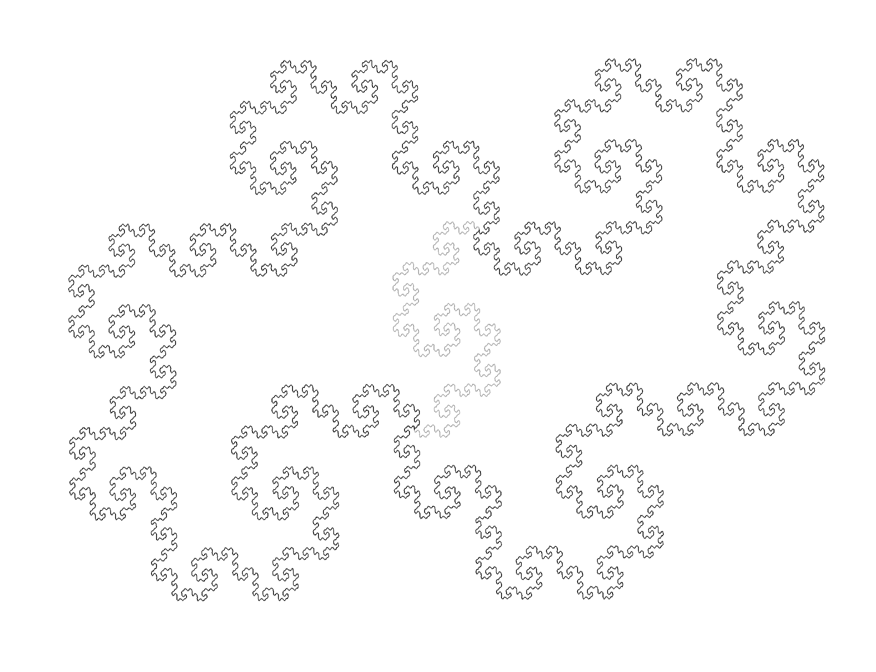}
  \end{center}
  \caption{The boundary of the twindragon fractal; \copyright~Wikipedia.}
  \label{Golden Dragon fractal}
\end{figure*}

In addition, the boundary in Figure \ref{fig:tornado_fractal} exhibits fractal structure, similar to that of the twindragon shown in Figure \ref{Golden Dragon fractal} as well as other ``dragon'' fractals. Interestingly, the fractal dimensions of the boundaries of the dragon fractals are in the range of the exponents in the power laws discussed in this work (e.g., $\sim1.52$ for the twindragon or $\sim1.62$ for the ``golden dragon'').

We recommend further exploration of the relationships among helicity, temperature (in the vortex gas sense), self-similarity, and the power laws proposed by Cai and Wurman. The resulting benefits of our understanding of helical atmospheric vortices could improve operational tornadic prediction. To the degree that the vorticity power law extends from observable to tornadic scales, it may also be possible to improve tornado detection and perhaps even estimate maximum tangential winds in tornadoes, as discussed in \citet{cai}. These hypotheses should be tested using real radar observations of tornadic and nontornadic supercells.

\begin{acknowledgment}
The four authors, B\v{e}l\'{\i}k, Dokken, Scholz, and Shvartsman were supported by National Science Foundation grant DMS-0802959. Funding for Potvin was provided by the NOAA/Office of Oceanic and Atmospheric Research under NOAA--University of Oklahoma Cooperative Agreement \#NA11OAR4320072, U.S.~Department of Commerce. Funding for Dahl and McGovern was provided by NSF/IIS grant 0746816.
\end{acknowledgment}

\ifthenelse{\boolean{dc}}
{}
{\clearpage}

\bibliographystyle{ametsoc} 
\bibliography{vortex_gas_selfsimilarity.01.15.2015}

\begin{thebibliography}{110}
\providecommand{\natexlab}[1]{#1}
\providecommand{\url}[1]{\texttt{#1}}
\providecommand{\urlprefix}{URL }
\expandafter\ifx\csname urlstyle\endcsname\relax
  \providecommand{\doi}[1]{doi:\discretionary{}{}{}#1}\else
  \providecommand{\doi}{doi:\discretionary{}{}{}\begingroup
  \urlstyle{rm}\Url}\fi
\providecommand{\eprint}[2][]{\url{#2}}

\bibitem[{Adlerman and Droegemeier(2000)}]{adlerman00}
Adlerman, E.~J. and K.~Droegemeier, 2000: A numerical simulation of cyclic
  tornadogenesis. \textit{20th Conference on Severe Local Storms}, American
  Meteorological Society, Orlando, FL, 591.

\bibitem[{Adlerman and Droegemeier(2002)}]{adlerman02}
Adlerman, E.~J. and K.~Droegemeier, 2002: The sensitivity of numerically
  simulated cyclic mesocyclogenesis to variations in environmental parameters.
  \textit{21th Conference on Severe Local Storms}, American Meteorological
  Society, San Antonio, TX.

\bibitem[{Andr\'{e} and Lesieur(1977)}]{andre77}
Andr\'{e}, J.~C. and M.~Lesieur, 1977: Influence of helicity on the evolution
  of isotropic turbulence at high {R}eynolds number. \textit{J. Fluid Mech.},
  \textbf{81~(1)}, 187--207.

\bibitem[{Arnold and Khesin(1999)}]{arnold99}
Arnold, V.~I. and B.~A. Khesin, 1999: \textit{Topological Methods in
  Hydrodynamics}, Applied Mathematical Sciences, Vol. 125. 2d ed.,
  Springer--Verlag, New York.

\bibitem[{Baker and Shelley(1990)}]{baker90}
Baker, G.~R. and M.~J. Shelley, 1990: On the connection between thin vortex
  layers and vortex sheet. \textit{J. Fluid Mech.}, \textbf{215}, 161--194.

\bibitem[{Barcilon(1967)}]{barcilon67}
Barcilon, A.~I., 1967: Vortex decay above a stationary boundary. \textit{J.
  Fluid Mech.}, \textbf{27~(1)}, 155--157.

\bibitem[{Barenblatt(1996)}]{barenblatt96}
Barenblatt, G.~I., 1996: \textit{Scaling, self-similarity, and intermediate
  asymptotics}. Cambridge University Press.

\bibitem[{Barenblatt(2003)}]{barenblatt03}
Barenblatt, G.~I., 2003: \textit{Scaling}. Cambridge University Press.

\bibitem[{Barnsley(1988)}]{barnsley88}
Barnsley, M., 1988: \textit{Fractals Everywhere}. Academic Press.

\bibitem[{Benjamin(1962)}]{benjamin62}
Benjamin, T.~B., 1962: Theory of the vortex breakdown phenomenon. \textit{J.
  Fluid Mech.}, \textbf{14}, 593--629.

\bibitem[{Bluestein(2013)}]{bluestein13}
Bluestein, H.~B., 2013: \textit{Severe Convective Storms and Tornadoes,
  Observations and Dynamics}. Springer-Praxis books in Environmental Sciences,
  Springer.

\bibitem[{Bluestein et~al.(2000)Bluestein, Gaddy, Dowell, Pazmany, Galloway,
  Mcintosh, and Stein}]{bluestein00a}
Bluestein, H.~B., S.~G. Gaddy, D.~C. Dowell, A.~L. Pazmany, J.~C. Galloway,
  R.~E. Mcintosh, and H.~Stein, 2000: Doppler radar observations of
  substorm-scale vortices in a supercell. \textit{Monthly Weather Review},
  \textbf{125}, 1046--1059.

\bibitem[{Bluestein et~al.(2003{\natexlab{a}})Bluestein, Lee, Bell, Weiss, and
  Pazmany}]{bluestein03b}
Bluestein, H.~B., W.~C. Lee, M.~Bell, C.~C. Weiss, and A.~L. Pazmany,
  2003{\natexlab{a}}: Mobile {D}oppler radar observations of a tornado in a
  supercell near {B}assett, {N}ebraska, on 5 {J}une 1999. {P}art {II}:
  Tornado-vortex. \textit{Mon. Wea. Rev.}, \textbf{131}, 2968--2984.

\bibitem[{Bluestein and Pazmany(2000)}]{bluestein00b}
Bluestein, H.~B. and A.~L. Pazmany, 2000: Observations of tornadoes and other
  convective phenomena with a mobile, 3-mm wavelength, {D}oppler radar: The
  spring 1999 field experiment. \textit{Bull. Amer. Met. Soc}, \textbf{81},
  2939--2951.

\bibitem[{Bluestein et~al.(2003{\natexlab{b}})Bluestein, Pazmany, and
  Weiss}]{bluestein03a}
Bluestein, H.~B., A.~L. Pazmany, and C.~C. Weiss, 2003{\natexlab{b}}: Mobile
  {D}oppler radar observations of a tornado in a supercell near {B}assett,
  {N}ebraska, on 5 {J}une 1999. {P}art {I}: {T}ornadogenesis. \textit{Mon. Wea.
  Rev.}, \textbf{131}, 2954--2967.

\bibitem[{Bohac(2013)}]{elrenovideo13}
Bohac, R., 2013: 2013 {EF5} {E}l {R}eno, {OK} tornado showing multiple
  funnels/subvortices. http://www.youtube.com/watch?v=C34EVyWRZbk.

\bibitem[{Brandes(1978)}]{brandes78}
Brandes, E.~A., 1978: Mesocyclone evolution and tornadogenesis: {S}ome
  observations. \textit{Mon. Wea. Rev.}, \textbf{106}, 995--1011.

\bibitem[{Bryan and Fritsch(2002)}]{bryanfritsch02}
Bryan, G.~H. and J.~M. Fritsch, 2002: A benchmark simulation for moist
  hydrodynamic numerical models. \textit{Mon. Wea. Rev.}, \textbf{130},
  2917--2928.

\bibitem[{Burggraf and Foster(1977)}]{Burggraf77}
Burggraf, O.~R. and M.~R. Foster, 1977: Continuation or breakdown in
  tornado-like vortices. \textit{J. Fluid Mech.}, \textbf{80~(4)}, 685--703.

\bibitem[{B\v{e}l\'{\i}k et~al.(2014)B\v{e}l\'{\i}k, Dokken, Scholz, and
  Shvartsman}]{belik14}
B\v{e}l\'{\i}k, P., D.~P. Dokken, K.~Scholz, and M.~M. Shvartsman, 2014:
  Fractal powers in {S}errin's swirling vortex model. \textit{Asymptot. Anal.},
  to appear.

\bibitem[{Cai(2005)}]{cai}
Cai, H., 2005: Comparison between tornadic and nontornadic mesocyclones using
  the vorticity (pseudovorticity) line technique. \textit{Mon. Wea. Rev.},
  \textbf{133~(9)}, 2535--2551.

\bibitem[{Chorin(1993)}]{chorin93}
Chorin, A.~J., 1993: Hairpin removal in vortex interactions {II}. \textit{J.
  Comp. Phys}, \textbf{107}, 1--9.

\bibitem[{Chorin(1994)}]{chorin}
Chorin, A.~J., 1994: \textit{Vorticity and turbulence}. Springer--Verlag, New
  York.

\bibitem[{Chorin and Akao(1991)}]{chorin-akao91}
Chorin, A.~J. and J.~Akao, 1991: Vortex equilibria in turbulence and quantum
  analogues. \textit{Physica D}, \textbf{52}, 403--414.

\bibitem[{Chorin and Bernard(1973)}]{chorin73}
Chorin, A.~J. and P.~Bernard, 1973: Discretization of a vortex sheet, with an
  example of roll-up. \textit{J. Comput. Phys.}, \textbf{13}, 423--429.

\bibitem[{Chorin and Marsden(1993)}]{chorin-marsden93}
Chorin, A.~J. and J.~E. Marsden, 1993: \textit{A Mathematical Introduction to
  Fluid Dynamics}. 3d ed., Springer--Verlag.

\bibitem[{Church et~al.(1977)Church, Snow, and Agee}]{church77}
Church, C.~R., J.~T. Snow, and E.~M. Agee, 1977: Tornado vortex simulation at
  {P}urdue {U}niversity. \textit{Bull. Amer. Meteor. Soc.}, \textbf{58~(9)},
  900--908.

\bibitem[{Cressman(1959)}]{cressman59}
Cressman, G.~P., 1959: An operational objective in analysis system.
  \textit{Mon. Wea. Rev.}, \textbf{87}, 367--374.

\bibitem[{Davies-Jones(1982)}]{davies-jones82}
Davies-Jones, R.~P., 1982: A new look at the vorticity equation with
  application to tornadogenesis. \textit{12th Conference on Severe Local
  Storms}, San Antonio, TX, Amer. Meteor. Soc., 249--252.

\bibitem[{Davies-Jones(1984)}]{davies-jones84}
Davies-Jones, R.~P., 1984: Streamwise vorticity: the origin of updraft rotation
  in supercell storms. \textit{J. Atmos. Sci.}, \textbf{41}, 2991--3006.

\bibitem[{Davies-Jones(1996)}]{davies-jones96}
Davies-Jones, R.~P., 1996: Formulas for the baroclinic and barotropic
  components of vorticity with application to vortex formation near the ground.
  \textit{Preprints, Seventh Conf. on Mesoscale Processes, Reading, United
  Kingdom}, Amer. Meteor. Soc., 14--16.

\bibitem[{Davies-Jones(2000)}]{davies-jones00}
Davies-Jones, R.~P., 2000: A {L}agrangian model for baroclinic genesis of
  mesoscale vortices. {P}art {I}: {T}heory. \textit{J. Atmos. Sci.},
  \textbf{57}, 715--736.

\bibitem[{Davies-Jones(2006{\natexlab{a}})}]{davies-jones06a}
Davies-Jones, R.~P., 2006{\natexlab{a}}: Integrals of the vorticity equation.
  {P}art {I}: {G}eneral three- and two-dimensional flows. \textit{J. Atmos.
  Sci.}, \textbf{63}, 598--610.

\bibitem[{Davies-Jones(2006{\natexlab{b}})}]{davies-jones06b}
Davies-Jones, R.~P., 2006{\natexlab{b}}: Integrals of the vorticity equation.
  {P}art {II}: {S}pecial two-dimensional flows. \textit{J. Atmos. Sci.},
  \textbf{63}, 611--616.

\bibitem[{Davies-Jones(2008)}]{davies-jones08}
Davies-Jones, R.~P., 2008: Can a descending rain curtain in a supercell
  instigate tornadogenesis barotropically? \textit{J. Atmos. Sci.},
  \textbf{65}, 2469--2497.

\bibitem[{Davies-Jones et~al.(1990)Davies-Jones, Burgess, and
  Foster}]{davis-jones90}
Davies-Jones, R.~P., D.~Burgess, and M.~Foster, 1990: Helicity as a tornado
  forecast parameter. \textit{16th Conference on Severe Local Storms},
  588--592.

\bibitem[{Deardorff(1980)}]{deardorff80}
Deardorff, J.~W., 1980: Stratocumulus-capped mixed layer derived from a
  three-dimensional model. \textit{Boundary-Layer Meteorology}, \textbf{18},
  495--527.

\bibitem[{Discovery.com/stormchasers(2010)}]{wadena10}
Discovery.com/stormchasers, 2010: Violent {M}innesota wedge tornado
  intercept!!! http://www.youtube.com/watch?v=AvD2nDyXSQo.

\bibitem[{Dowell and Bluestein(1997)}]{dowell97}
Dowell, D.~C. and H.~B. Bluestein, 1997: The {A}rcadia, {O}klahoma storm of 17
  {M}ay 1981: {A}nalysis of a supercell during tornadogenesis. \textit{Mon.
  Wea. Rev.}, \textbf{125}, 2562--2582.

\bibitem[{Dowell and Bluestein(2002{\natexlab{a}})}]{dowell02a}
Dowell, D.~C. and H.~B. Bluestein, 2002{\natexlab{a}}: The 8 {J}une 1995
  {M}c{L}ean, {T}exas, storm. {P}art {I}: {O}bservations of cyclic
  tornadogenesis. \textit{Mon. Wea. Rev.}, \textbf{125}, 2626--2648.

\bibitem[{Dowell and Bluestein(2002{\natexlab{b}})}]{dowell02b}
Dowell, D.~C. and H.~B. Bluestein, 2002{\natexlab{b}}: The 8 {J}une 1995
  {M}c{L}ean, {T}exas, storm. {P}art {II}: {C}yclic tornado formation,
  maintenance, and dissipation. \textit{Mon. Wea. Rev.}, \textbf{130},
  2649--2670.

\bibitem[{Droegemeier et~al.(1993)Droegemeier, Lazarus, and
  Davies-Jones}]{droegemeier93}
Droegemeier, K.~K., S.~M. Lazarus, and R.~P. Davies-Jones, 1993: The influence
  of helicity on numerically simulated convective storms. \textit{Mon. Wea.
  Rev.}, \textbf{121}, 2005--2029.

\bibitem[{Fiedler(1994)}]{fiedler94}
Fiedler, B.~H., 1994: The thermodynamic speed limit and its violation in
  axisymmetric numerical simulations of tornado-like vortices. \textit{Atmos.
  Ocean}, \textbf{32~(2)}, 335--359.

\bibitem[{Fiedler(1997)}]{fiedler97}
Fiedler, B.~H., 1997: Compressibility and windspeed limits in tornadoes.
  \textit{Atmos.--Ocean}, \textbf{35}, 93--107.

\bibitem[{Fiedler and Rotunno(1986)}]{fiedlerrotunno86}
Fiedler, B.~H. and R.~Rotunno, 1986: A theory for the maximum windspeed in
  tornado-like vortices. \textit{J. Atmos. Sci.}, \textbf{43~(21)}, 2328--2440.

\bibitem[{Flandoli and Gubinelli(2002)}]{flandoli02}
Flandoli, F. and M.~Gubinelli, 2002: The {G}ibbs ensemble of a vortex filament.
  \textit{Prob. Th. Rel. Fields}, \textbf{112}.

\bibitem[{Frehlich and Sharman(2008)}]{frehlich2008}
Frehlich, R. and R.~Sharman, 2008: The use of structure functions and spectra
  from numerical model output to determine effective model resolution.
  \textit{Mon. Wea. Rev.}, \textbf{136}, 1537--1553.

\bibitem[{Fr\"{o}hlich and Ruelle(1982)}]{frohlich-ruelle82}
Fr\"{o}hlich, J. and D.~Ruelle, 1982: Statistical mechanics of vortices in an
  inviscid two-dimensional fluid. \textit{Commun. Math. Phys.},
  \textbf{87~(1)}, 1--36.

\bibitem[{Fujita(1973)}]{fujita73}
Fujita, T.~T., 1973: Proposed mechanism of tornado formation from rotating
  thunderstorm. \textit{8th Conference on Severe Local Storms}, Denver, CO,
  Amer. Meteor. Soc., 191--196.

\bibitem[{Fujita(1975)}]{fujita75}
Fujita, T.~T., 1975: New evidence from the {A}pril 3--4, 1974 tornadoes.
  \textit{9th Conference on Severe Local Storms}, Norman, OK, Amer. Meteor.
  Soc., 248--255.

\bibitem[{Fujita(1981)}]{fujita81}
Fujita, T.~T., 1981: Tornadoes and downbursts in the context of generalized
  planetary scales. \textit{J. Atmos. Sci.}, \textbf{38~(8)}, 1511--1534.

\bibitem[{Garrod(1995)}]{garrod95}
Garrod, C., 1995: \textit{Statistical Mechanics and Thermodynamics}. Oxford U.
  Press.

\bibitem[{Grazulis(1997)}]{grazulis97}
Grazulis, T.~P., 1997: \textit{Significant Tornadoes Update 1992--1995}.
  Environmental Films, Copyright January 1997, P.O. Box 302, St. Johnsbury, VT
  05819.

\bibitem[{Helmholtz(1858)}]{helmholtz58}
Helmholtz, 1858: \"{U}ber {I}ntegrale der hydrodynamischen {G}leichungen welche
  den {W}irbelbewegungen entsprechen. \textit{Crelle}, \textbf{55}, 25--55.

\bibitem[{Katok and Hasselblatt(1995)}]{katok-hasselblatt95}
Katok, A. and B.~Hasselblatt, 1995: \textit{Introduction to the Modern Theory
  of Dynamical Systems}, Encyclopedia of Mathematics and its Applications,
  Vol.~54. Cambridge University Press.

\bibitem[{Klemp(1987)}]{klemp87}
Klemp, J.~B., 1987: Dynamics of tornadic thunderstorms. \textit{Ann. Rev. Fluid
  Mech.}, \textbf{19}, 369--402.

\bibitem[{Krasny(1993)}]{krasny93}
Krasny, R., 1993: Vortex sheet roll-up. \textit{RIMS Workshop on Unstable and
  Turbulent Motion of Fluid, Kyoto, Japan}, S.~Kida, Ed., World Scientific,
  43--49.

\bibitem[{Kuznetsov(2011)}]{kuznetsov11}
Kuznetsov, S.~P., 2011: Dynamical chaos and uniformly hyperbolic attractors:
  from mathematics to physics. \textit{Physics--Uspekhi}, \textbf{54~(2)},
  119--144.

\bibitem[{Lee and Wilhelmson(1997{\natexlab{a}})}]{lee-wilhelmson97a}
Lee, B.~D. and R.~B. Wilhelmson, 1997{\natexlab{a}}: The numerical simulation
  of non-supercell tornadogenesis: {P}art {I}: {I}nitiation and evolution of
  pretornadic misocyclone circulations along a dry outflow boundary. \textit{J.
  Atmos. Sci.}, \textbf{54}, 32--60.

\bibitem[{Lee and Wilhelmson(1997{\natexlab{b}})}]{lee-wilhelmson97b}
Lee, B.~D. and R.~B. Wilhelmson, 1997{\natexlab{b}}: The numerical simulation
  of non-supercell tornadogenesis, {P}art {II}: {E}volution of a family of
  tornadoes along a weak outflow boundary. \textit{J. Atmos. Sci.},
  \textbf{54}, 2387--2415.

\bibitem[{Lee and Wilhelmson(2000)}]{lee-wilhelmson00}
Lee, B.~D. and R.~B. Wilhelmson, 2000: The numerical simulation of
  non-supercell tornadogenesis, {P}art {III}: {P}arameter tests investigating
  the role of {CAPE}, vortex sheet strength and boundary layer vertical shear.
  \textit{J. Atmos. Sci.}, \textbf{57}, 2246--2261.

\bibitem[{Levich and Tzvetkov(1985)}]{levich85}
Levich, E. and E.~Tzvetkov, 1985: Helical inverse cascade in three-dimensional
  turbulence as a fundamental dominant mechanism in mesoscale atmospheric
  phenomena. \textit{Physics reports}, \textbf{128~(1)}, 1--37.

\bibitem[{Lewellen and Lewellen(2007{\natexlab{a}})}]{lewellens07}
Lewellen, D.~C. and W.~S. Lewellen, 2007{\natexlab{a}}: Near-surface
  intensification of tornado vortices. \textit{J. Atmos. Sci.},
  \textbf{64~(7)}, 2176--2194.

\bibitem[{Lewellen and Lewellen(2007{\natexlab{b}})}]{lewellens07a}
Lewellen, D.~C. and W.~S. Lewellen, 2007{\natexlab{b}}: Near-surface vortex
  intensification through corner flow collapse. \textit{J. Atmos. Sci.},
  \textbf{64~(7)}, 2195--2209.

\bibitem[{Lewellen et~al.(2000)Lewellen, Lewellen, and Xia}]{lewellensxia00}
Lewellen, D.~C., W.~S. Lewellen, and J.~Xia, 2000: The influence of a local
  swirl ratio on tornado intensification near the surface. \textit{J. Atmos.
  Sci.}, \textbf{57~(4)}, 527--544.

\bibitem[{Lilly(1969)}]{lilly69}
Lilly, D.~K., 1969: Tornado dynamics, {NCAR} {M}anuscript 69-117.

\bibitem[{Lilly(1976)}]{lilly76}
Lilly, D.~K., 1976: Sources of rotation and energy in the tornado.
  \textit{Proceedings of the Symposium on Tornadoes, Assessment of Knowledge
  and Implications for Man}, Lubbock, TX, Texas Tech University, 145--150.

\bibitem[{Lilly(1983)}]{lilly83}
Lilly, D.~K., 1983: Dynamics of rotating thunderstorms. \textit{Mesoscale
  Meteorology--Theories, Observations, and Models}, D.~K. Lilly and
  E.~T.~Gal-Chen, Eds., D. Reidel, Dordrecht, Vol. 114, 531--544.

\bibitem[{Lilly(1986{\natexlab{a}})}]{lilly86a}
Lilly, D.~K., 1986{\natexlab{a}}: The structure, energetics and propagation of
  rotating convective storms. {P}art {I}: {E}nergy exchange with the mean flow.
  \textit{J. Atmos. Sci.}, \textbf{43}, 113--125.

\bibitem[{Lilly(1986{\natexlab{b}})}]{lilly86b}
Lilly, D.~K., 1986{\natexlab{b}}: The structure, energetics and propagation of
  rotating convective storms. {P}art {II}: {H}elicity and storm stablization.
  \textit{J. Atmos. Sci.}, \textbf{43}, 126--140.

\bibitem[{Lim and Nebus(2007)}]{lim-nebus07}
Lim, C. and J.~Nebus, 2007: \textit{Vorticity, Statistical Mechanics, and
  {M}onte {C}arlo Simulation}. Springer--Verlag.

\bibitem[{Lund and Snow(1993)}]{lund-snow93}
Lund, D.~E. and J.~T. Snow, 1993: Laser {D}oppler velocimeter mesaurements in
  tornado-like vortices. \textit{The Tornado: Its Structure, Dynamics,
  Prediction, and Hazards.}, American Geophysiscal Union, Vol. Monograph 79,
  297--306.

\bibitem[{Majda and Bertozzi(2001)}]{majda01}
Majda, A.~J. and A.~Bertozzi, 2001: \textit{Vorticity and Incompressible
  Flows}. 1st ed., Cambridge Texts in Applied Mathematics, Cambridge University
  Press.

\bibitem[{Mandelbrot(1983)}]{mandelbrot83}
Mandelbrot, B.~B., 1983: \textit{Tha fractal geometry of nature}. W.H. Freeman
  and Company.

\bibitem[{Marchioro and Pulvirenti(1994)}]{marchioro94}
Marchioro, C. and M.~Pulvirenti, 1994: \textit{Mathematical Theory of
  Incompressible Nonviscous Fluids}, Applied Mathematical Sciences, Vol.~96.
  Springer--Verlag.

\bibitem[{Markowski et~al.(2002)Markowski, Straka, and
  Rasmussen}]{markowski02a}
Markowski, P.~A., J.~M. Straka, and E.~N. Rasmussen, 2002: Direct surface
  thermodynamic observations within the rear-flank downdrafts of nontornadic
  and tornadic supercells. \textit{Mon. Wea. Rev.}, \textbf{130}, 1692--1721.

\bibitem[{Markowski et~al.(2003)Markowski, Straka, and Rasmussen}]{markowski03}
Markowski, P.~M., J.~M. Straka, and E.~N. Rasmussen, 2003: Tornadogenesis
  resulting from transport of circulation by a downdraft: {I}dealized numerical
  simulations. \textit{J. Atmos. Sci.}, \textbf{60}, 795--823.

\bibitem[{Markowski et~al.(2008)Markowski, Straka, Rasmussen, Davies-Jones,
  Richardson, and Trapp}]{markowski08}
Markowski, P.~M., J.~M. Straka, E.~N. Rasmussen, R.~P. Davies-Jones,
  Y.~Richardson, and J.~Trapp, 2008: Vortex lines within low-level mesocyclones
  obtained from pseudo-dual-{D}oppler radar observations. \textit{Mon. Wea.
  Rev.}, \textbf{136}, 3513--3535.

\bibitem[{Miller(1967)}]{miller67}
Miller, B.~I., 1967: Characteristics of hurricanes. \textit{Science},
  \textbf{157~(3795)}, 1389--1399.

\bibitem[{Miller et~al.(1992)Miller, Weichman, and Cross}]{miller92}
Miller, J., P.~B. Weichman, and M.~C. Cross, 1992: Statistical mechanics,
  {E}uler's equation, and {J}upiter's {R}ed {S}pot. \textit{Physical Review A},
  \textbf{45}, 2328--2359.

\bibitem[{Moffat(1969)}]{moffat69}
Moffat, H.~K., 1969: The degree of knottedness of tangled vortex lines.
  \textit{J. Fluid Mech.}, \textbf{35~(1)}, 117--129.

\bibitem[{Moffat and Tsinober(1992)}]{moffat92}
Moffat, H.~K. and A.~Tsinober, 1992: Helicity in laminar and turbulent flow.
  \textit{Annu. Rev. Fluid Mech.}, \textbf{24}, 281--312.

\bibitem[{Morrison et~al.(2005)Morrison, Curry, and Khvorostyanov}]{morrison05}
Morrison, H., J.~A. Curry, and V.~I. Khvorostyanov, 2005: A new double-moment
  microphysics parametrization for application in cloud and climate models.
  {P}art {I}: Description. \textit{J. Atmos. Sci}, \textbf{62}, 1665--1677.

\bibitem[{Newton(2001)}]{newton01}
Newton, P.~K., 2001: \textit{The {N}-vortex problem. {A}nalytical Techniques}.
  Springer--Verlag, New York.

\bibitem[{No\-va(2004)}]{nova04}
No\-va, 2004: Hunt for the super twister.
  http://www.pbs.org/wgbh/nova/earth/hunt-for-the-supertwister.html.

\bibitem[{Nolan(2012)}]{nolan12}
Nolan, D.~S., 2012: Three-dimensional instabilities in tornado-like vortices
  with secondary circulations. \textit{J. Fluid Mech.}, \textbf{711}, 61--100.

\bibitem[{Onsager(1949)}]{onsager49}
Onsager, L., 1949: Statistical hydrodynamics. \textit{Il Nuovo Cimento},
  \textbf{6}, 279--287.

\bibitem[{Pesin and Climenhaga(2009)}]{pesin09}
Pesin, Y. and V.~Climenhaga, 2009: \textit{Lectures on {F}ractal {G}eometry and
  {D}ynamical {S}ystems}, Student Mathematical Library, Vol.~52. AMS,
  Providence, RI.

\bibitem[{Potvin(2013)}]{potvin13}
Potvin, C.~K., 2013: A variational method for detecting and characterizing
  intense vortices in {C}artesian wind fields. \textit{Mon. Wea. Rev.},
  \textbf{141}, 3102--3115.

\bibitem[{Pouquet and Mininni(2010)}]{pouquet10}
Pouquet, A. and P.~D. Mininni, 2010: The interplay between helicity and
  rotation in turbulence: implications for scaling laws and small-scale
  dynamics. \textit{Phil. Trans. R. Soc. A}, \textbf{368}, 1635--1662.

\bibitem[{Roberts and Wilson(1995)}]{roberts95}
Roberts, R.~D. and J.~W. Wilson, 1995: The genesis of three nonsupercell
  tornadoes observed with dual-{D}oppler radar. \textit{Mon. Wea. Rev.},
  \textbf{123}, 3408--3436.

\bibitem[{Rotunno(1984)}]{rotunno84}
Rotunno, R., 1984: An investigation of a three-dimensional asymmetric vortex.
  \textit{J. Atmos. Sci.}, \textbf{41}, 283--298.

\bibitem[{Schwartz(1976)}]{lschwartz76}
Schwartz, L., 1976: Lectures on disintegration of measures. \textit{Tata
  Institute of Fundamental Research}.

\bibitem[{Serrin(1972)}]{serrin}
Serrin, J., 1972: The swirling vortex. \textit{Phil. Trans. Roy. Soc. London,
  Series A, Math \& Phys. Sci.}, \textbf{271~(1214)}, 325--360.

\bibitem[{Snow(1978)}]{snow78}
Snow, J.~T., 1978: On inertial instability as related to the multiple vortex
  phenomena. \textit{J. Atmos. Sci.}, \textbf{35}, 1660--1677.

\bibitem[{Straka et~al.(2007)Straka, Rasmussen, Davies-Jones, and
  Markowski}]{straka07}
Straka, J.~M., E.~N. Rasmussen, R.~P. Davies-Jones, and P.~M. Markowski, 2007:
  An observational and idealized numerical examination of low-level
  counter-rotating vortices toward the rear flank of supercells.
  \textit{Electronic J. Severe Storms Meteor.}, \textbf{2~(8)}, 1--22.

\bibitem[{Thomson(1869)}]{kelvin69}
Thomson, L. K.~W., 1869: On vortex motion. \textit{Trans. Roy. Soc. Edinburgh},
  \textbf{25}, 217--260.

\bibitem[{Trapp(1999)}]{trapp99}
Trapp, R.~J., 1999: Observations of nontornadic low-level mesocyclones and
  attendant tornadogenesis failure during {VORTEX} 94. \textit{Mon. Wea. Rev.},
  \textbf{124}, 384--407.

\bibitem[{Trevorrow et~al.(2012)Trevorrow, Hildner, Tripoli, and
  B{\"u}ker}]{trevorrow12}
Trevorrow, S.~T., R.~D. Hildner, G.~J. Tripoli, and M.~L. B{\"u}ker, 2012:
  Numerical study of low-level vorticity prior to tornadogenesis. \textit{26th
  Conference on Severe Local Storms}, American Meteorological Society,
  Nashville, TN, Amer. Meteor. Soc.

\bibitem[{Wakimoto and Atkins(1996)}]{wakimoto96}
Wakimoto, R.~M. and N.~T. Atkins, 1996: Observations on the origins of
  rotation: {T}he {N}ewcastle tornado during {VORTEX} 94. \textit{Mon. Wea.
  Rev.}, \textbf{124}, 384--407.

\bibitem[{Wakimoto and Wilson(1989)}]{wakimoto89}
Wakimoto, R.~M. and J.~W. Wilson, 1989: Non-supercell tornadoes. \textit{Mon.
  Wea. Rev.}, \textbf{117}, 1113--1140.

\bibitem[{Waleffe(1992)}]{waleffe92}
Waleffe, F., 1992: The nature of triad interactions in homogeneous turbulence.
  \textit{Phys. Fluids A}, \textbf{4~(2)}, 350--363.

\bibitem[{Ward(1972)}]{ward72}
Ward, N.~B., 1972: The exploration of certain features of tornado dynamics
  using a laboratory model. \textit{J. Atmos. Sci.}, \textbf{29~(6)},
  1194--1204.

\bibitem[{Wilson(1986)}]{wilson86}
Wilson, J.~W., 1986: Tornadogenesis by nonprecipitation induced wind shear
  lines. \textit{Mon. Wea. Rev.}, \textbf{114}, 270--284.

\bibitem[{Wurman(2002)}]{wurman02}
Wurman, J., 2002: The multiple-vortex structure of a tornado. \textit{Weather
  and Forecasting}, \textbf{17}, 473--505.

\bibitem[{Wurman and Alexander(2005)}]{wurman05}
Wurman, J. and C.~R. Alexander, 2005: The 30 {M}ay 1998 {S}pencer, {S}outh
  {D}akota, storm. {P}art {II}: {C}omparison of observed damage and
  radar-derived winds in the tornadoes. \textit{Mon. Wea. Rev.},
  \textbf{133~(1)}, 97--119.

\bibitem[{Wurman and Gill(2000)}]{wurman00}
Wurman, J. and S.~Gill, 2000: Fine-scale radar observations of the {D}immitt,
  {T}exas (2 {J}une 1995) tornado. \textit{Mon. Wea. Rev.}, \textbf{128},
  2135--2164.

\bibitem[{Wurman et~al.(2013)Wurman, Kosiba, and Robinson}]{wurman2013}
Wurman, J., K.~Kosiba, and P.~Robinson, 2013: In situ, doppler radar, and video
  observations of the interior structure of a tornado and the wind-damage
  relationship. \textit{Bull. Amer. Meteor. Soc.}, \textbf{94~(6)}, 835--846.

\bibitem[{Xia et~al.(2003)Xia, Lewellen, and Lewellen}]{xialewellens03}
Xia, J., D.~C. Lewellen, and W.~S. Lewellen, 2003: Influence of {M}ach number
  on tornado corner flow dynamics. \textit{J. Atmos. Sci.}, \textbf{60~(22)},
  2820--2825.

\bibitem[{Yokoi and Yoshizawa(1993)}]{yokoi93}
Yokoi, N. and A.~Yoshizawa, 1993: Statistical analysis of the effects of
  helicity in inhomogeneous turbulence. \textit{Phys. Fluids A},
  \textbf{5~(2)}, 464--477.

\end{thebibliography}


\end{document}